\documentclass[11pt,a4paper,twoside]{article}
\usepackage[T1]{fontenc}
\usepackage{lmodern}
\usepackage[a4paper,left=25mm,right=25mm,top=24mm,bottom=25mm,headheight=14pt,headsep=8mm]{geometry}
\usepackage{amsmath,amssymb,amsthm,mathtools}
\usepackage{microtype}
\usepackage{booktabs,longtable,array,enumitem}
\usepackage{xcolor}
\usepackage{fancyhdr}
\fancypagestyle{plain}{%
  \fancyhf{}
  \fancyfoot[C]{\scriptsize\thepage}
  
  }
\renewcommand{\footnoterule}{%
  \kern-3pt\hrule width 5pc\kern2.6pt}
\usepackage{hyperref}
\hypersetup{colorlinks=true,linkcolor=blue!45!black,citecolor=blue!45!black,
  urlcolor=blue!45!black,pdftitle={Compactness from fractional momentum bounds and weak solutions of compressible Navier-Stokes},
  pdfauthor={Yanlin Liu and Ping Zhang},bookmarksnumbered=true,bookmarksdepth=3}
\numberwithin{equation}{section}
\newtheorem{theorem}{Theorem}[section]
\newtheorem{lemma}[theorem]{Lemma}
\newtheorem{proposition}[theorem]{Proposition}
\newtheorem{corollary}[theorem]{Corollary}
\theoremstyle{definition}
\newtheorem{definition}{Definition}[section]
\theoremstyle{remark}
\newtheorem{remark}{Remark}[section]
\newcommand{\R}{\mathbb R}

\newcommand{\dd}{\,\mathrm d}
\newcommand{\diver}{\operatorname{div}}

\newcommand{\tr}{\operatorname{tr}}

\newcommand{\supp}{\operatorname{supp}}

\newcommand{\esssup}{\operatorname*{ess\,sup}}
\setlist{topsep=4pt,itemsep=2pt,parsep=0pt}
\newcommand{\Div}{\operatorname{Div}}

\allowdisplaybreaks[2]
\title{\MakeUppercase{Compactness from fractional momentum bounds}\\
\MakeUppercase{and weak solutions of the compressible}\\
\MakeUppercase{Navier--Stokes equations}}

\makeatletter
\gdef\@author{}
\gdef\ns@authorblocks{}
\renewcommand{\author}[2][]{%
  \gdef\ns@currentauthor{#1}%
  \ifx\@author\@empty
    \gdef\@author{#2}%
  \else
    \g@addto@macro\@author{ and #2}%
  \fi
  \g@addto@macro\ns@authorblocks{\ns@printauthor{#1}}}
\newcommand{\address}[2][]{%
  \expandafter\gdef\csname ns@address@#1\endcsname{#2}}
\newcommand{\email}[1]{%
  \expandafter\gdef\csname ns@email@\ns@currentauthor\endcsname{#1}}
\newcommand{\ns@printauthor}[1]{%
  \par\addvspace{6pt}%
  \indent\textup{(#1)} {\scshape\csname ns@address@#1\endcsname}\par
  \indent\textit{Email address: }\texttt{\csname ns@email@#1\endcsname}\par}
\renewcommand{\maketitle}{%
  \par\begingroup
    \thispagestyle{plain}%
    \vspace*{12pt}%
    {\centering
      {\normalfont\normalsize\bfseries\@title\par}%
      \vspace{18pt}%
      {\normalfont\footnotesize\MakeUppercase{\@author}\par}%
      \vspace{20pt}}%
    \renewcommand{\thefootnote}{}%
    \begin{NoHyper}%
      \footnotetext{\textit{Date: }\@date.}%
    \end{NoHyper}%
  \endgroup
  \setcounter{footnote}{0}%
  \par}
\AtEndDocument{%
  \par\addvspace{12pt}%
  \noindent\begin{minipage}{\textwidth}
    \normalfont\footnotesize
    \setlength{\parindent}{12pt}%
    \setlength{\parskip}{0pt}%
    \ns@authorblocks
  \end{minipage}\par}
\makeatother

\author[Y. Liu]{Yanlin Liu}
\address[Y. Liu]{School of Mathematical Sciences,
Laboratory of Mathematics and Complex Systems,
MOE, Beijing Normal University, 100875 Beijing, China.}
\email{liuyanlin@bnu.edu.cn}

\author[P. Zhang]{Ping Zhang}
\address[P. Zhang]{Academy of Mathematics $\&$ Systems Science
	and Hua Loo-Keng Center for Mathematical Sciences, Chinese Academy of
	Sciences, Beijing 100190, CHINA, and School of Mathematical Sciences,
	University of Chinese Academy of Sciences, Beijing 100049, China.} \email{zp@amss.ac.cn}

\date{November 19, 2026}
\begin{document}
\maketitle
\begin{abstract}
We consider the barotropic Navier--Stokes equations on $\R^3$ with
constant viscosities, finite mass and finite energy. For every $\gamma>1$,
a local bound on the momentum $m=\rho u$ in $L^2_tH^{-1+\sigma}_x$,
$\sigma>0$, gives strong local convergence of the densities under the hypotheses of
Theorem~\ref{cp:compactness}. For $\gamma>9/8$ we prove, uniformly along
the approximations,
\[
 \|I_{1-\alpha}m\|_{L^2(I\times\R^3)}\le C_{I,\alpha}
 \qquad\text{for some }\alpha>0,\quad I\Subset(0,\infty).
\]
Here $I_s$ is the spatial Riesz potential. Together with tightness of
the energy at spatial infinity, this gives a
global finite-energy renormalized weak solution.
\end{abstract}

\section{Introduction}\label{sec:main}

Let $\rho\ge0$, $u\in\R^3$, and $m=\rho u$. On $(0,\infty)\times\R^3$,
the barotropic Navier--Stokes equations read
\begin{equation}\label{eq:NS}
 \begin{aligned}
  \partial_t\rho+\diver m&=0,\\
  \partial_t m+\diver(m\otimes u)+\nabla p(\rho)&=\diver S(u),
  &p(\rho)&=a\rho^\gamma,\quad a>0,\quad\gamma>1.
 \end{aligned}
\end{equation}
The viscous stress is
\begin{equation}\label{eq:stress}
 S(u)=\mu(\nabla u+\nabla u^{\top})
          +\lambda\,\diver u\,\mathrm{Id}_3,
 \qquad \mu>0,\quad \lambda+\tfrac23\mu\ge0 .
\end{equation}
The initial data have finite mass and energy:
\begin{equation}\label{eq:initial-data}
 \begin{gathered}
 \rho_0\ge0,\qquad \rho_0\in L^1(\R^3)\cap L^\gamma(\R^3),\qquad
 m_0=0\quad\text{a.e. on }\{\rho_0=0\},\\
 E_0=\int_{\R^3}\left(\frac{|m_0|^2}{2\rho_0}
             +\frac{a}{\gamma-1}\rho_0^\gamma\right)\dd x<\infty .
 \end{gathered}
\end{equation}
We set $|m|^2/\rho=0$ at $(0,0)$ and $+\infty$ if $\rho=0$, $m\ne0$.
Then $m_0\in L^1\cap L^{2\gamma/(\gamma+1)}$. No symmetry, smallness or
moment condition is imposed.

For general data in three dimensions, see Lions \cite{Lions1998}
($\gamma\ge9/5$) and Feireisl, Novotn\'y and Petzeltov\'a
\cite{FNP2001} ($\gamma>3/2$ on bounded regular domains).
In a recent preprint, Hong, Peng and Zhu \cite{HPZ2026}
report existence for $\gamma>4/3$ on the three-dimensional torus.
Their Remark~1.5 also states extensions to $\R^3$ and smooth bounded
domains. The present paper treats $\R^3$.
We use the effective viscous flux and density truncations as in
\cite{Feireisl2001,FNP2001}. For spherically symmetric data and
$\gamma>1$, see Jiang and Zhang \cite{JiangZhang2001}; for the axisymmetric
case, see \cite{JiangZhang2003}. Nonmonotone pressures and anisotropic
stresses are treated in \cite{BreschJabin2018}.

Unless otherwise stated, spatial norms are over $\R^3$.
The spaces $\dot H^s$ and $H^s$ are the homogeneous and inhomogeneous
Sobolev spaces. We write
\[
 \begin{aligned}
 E(t)&=\int_{\R^3}\left(\frac{|m|^2}{2\rho}
               +\frac{a}{\gamma-1}\rho^\gamma\right)(t)\dd x,\\
 \mathcal D(u)&=\int_{\R^3}
     \bigl(\mu|\nabla u|^2+(\mu+\lambda)|\diver u|^2\bigr)\dd x .
 \end{aligned}
\]
\begin{definition}[Weak solution and energy inequalities]\label{def:solution}
A \emph{bounded-energy renormalized weak solution} satisfies, for each $0<T<\infty$,
\begin{equation}\label{eq:energy-spaces}
 \begin{gathered}
 \rho\ge0,\qquad \rho\in L^\infty(0,T;L^1\cap L^\gamma),\\
 \sqrt\rho\,u\in L^\infty(0,T;L^2),\qquad
 u\in L^2(0,T;\dot H^1).
 \end{gathered}
\end{equation}
We take the $L^6$ representative of $\dot H^1(\R^3)$.
For smooth scalar tests $\varphi$ and vector tests $\Phi$ compactly
supported in $[0,\infty)\times\R^3$,
\begin{equation}\label{eq:weak-mass}
 \int_0^\infty\!\!\int_{\R^3}(\rho\partial_t\varphi+m\cdot\nabla\varphi)
       \dd x\dd t=-\int_{\R^3}\rho_0\varphi(0)\dd x ,
\end{equation}
\begin{equation}\label{eq:weak-momentum}
 \begin{aligned}
 \int_0^\infty\!\!\int_{\R^3}
  \bigl[m\cdot\partial_t\Phi+(m\otimes u):\nabla\Phi
      +p(\rho)\diver\Phi-S(u):\nabla\Phi\bigr]\dd x\dd t
  &=-\int_{\R^3}m_0\cdot\Phi(0)\dd x .
 \end{aligned}
\end{equation}
For $b\in C^1([0,\infty))$ with $b'(z)=0$ for all sufficiently large $z$,
\begin{equation}\label{eq:renormalized}
 \begin{aligned}
 \int_0^\infty\!\!\int_{\R^3}
 \bigl[b(\rho)\cdot(\partial_t\varphi+u\cdot\nabla\varphi)
     +(b(\rho)-\rho b'(\rho))\varphi\,\diver u\bigr]\dd x\dd t\\
 =-\int_{\R^3}b(\rho_0)\varphi(0)\dd x .
 \end{aligned}
\end{equation}
Finally,
\begin{equation}\label{eq:energy-zero}
 E(t)+\int_0^t\mathcal D(u(s))\dd s\le E_0
 \qquad\text{for a.e. }t>0 .
\end{equation}
A \emph{finite-energy renormalized weak solution} satisfies in addition
\begin{equation}\label{r4:energy-definition}
 -\int_0^\infty E(t)\varphi'(t)\dd t+
 \int_0^\infty\mathcal D(u(t))\varphi(t)\dd t\le0
 \qquad (0\le\varphi\in C_c^\infty(0,\infty)).
\end{equation}
\end{definition}
Inequality \eqref{r4:energy-definition} gives the energy inequality between
almost every pair of positive times. On $\R^3$ its passage to the limit
also requires tightness of the energy at infinity \cite{DellOroFeireisl2015}.

\paragraph{Notation.}
We always denote $C$ or $c>0$ to be an absolute constant
which may vary from line to line.
Parameter dependence is indicated by subscripts; any additional fixed data
included in a subscripted constant are specified locally. Estimates stated as uniform, and the constants in $\lesssim$, are independent
of the approximation indices and of the indicated truncation and
regularization parameters.
$I\Subset J$ means that $\overline I$ is a compact subset of $J$,
and $C_{\rm weak}$ denotes continuity in the weak topology.
The spatial Riesz potential $I_s$ acts componentwise:
\[
 I_s f(x)=c_s\int_{\R^3}\frac{f(y)}{|x-y|^{3-s}}\dd y,\qquad
 \widehat{I_s f}(\xi)=|\xi|^{-s}\widehat f(\xi),\qquad 0<s<3,
\]
where $c_s$ is the normalization constant corresponding to the stated Fourier multiplier.

\subsection{Main theorem: the range \texorpdfstring{$\gamma>9/8$}{gamma > 9/8}}
The main result of this paper is Theorem~\ref{thm:nine-eighths}. Its proof separates
two tasks.

\paragraph{Compactness of the nonlinear terms.}
At low pressure exponents, the natural energy bounds do not by
themselves identify the weak limits of the pressure and convective
flux. Theorem~\ref{cp:compactness} achieves this
for every $\gamma>1$, provided a positive fractional bound on the
localized momentum is available. Bounded density truncations first
identify the momentum; its additional compactness then permits the
identification of the convective flux and the effective viscous flux.

\paragraph{A uniform estimate from the equations.}
The momentum bound must itself be derived from the equations, with
constants that remain uniform as the artificial pressure vanishes.
Section~\ref{ao:section} obtains it for $\gamma>9/8$. The analytic
input is the projected determinant estimate of
Proposition~\ref{rt:optimized-CI}, which accommodates the divergence
of an $L^2$ viscous stress. A cover constructed from the momentum
equation then converts this estimate into fractional regularity.

We first state the compactness criterion. In addition to the equations
and energy bounds, it requires control of the localized momentum in
$L^2_tH^{-1+\sigma}_x$ for some $\sigma>0$.

\begin{theorem}[Compactness of approximations]\label{cp:compactness}
Let $\gamma>1$, $T_j\to\infty$, and let $\Omega_j=\R^3$ or let
$\Omega_j$ be balls exhausting $\R^3$. For $(\rho_j,u_j)$ on
$(0,T_j)\times\Omega_j$, set
\[
 \rho_j\ge0,\qquad m_j=\rho_j u_j,\qquad S_j=S(u_j).
\]
On balls, impose $u_j|_{\partial\Omega_j}=0$ and extend $\rho_j,u_j$,
as well as the initial density and momentum, by zero. For whole-space
integrals, form the pressure fields $P_j,\pi_j(\rho_j)$ and the
internal-energy densities $\Theta_j(\rho_j),\Theta_j(\rho_{0j})$ on
$\Omega_j$ and then extend them by zero. In the renormalized equation,
form $b(\rho_j)$ and its initial value $b(\rho_{0j})$ from the extended
densities, so their exterior value is $b(0)$. In either case $u_j$ denotes the
$L^6(\R^3)$ representative in $L^2(0,T;\dot H^1(\R^3))$ for
each $0<T<T_j$. In particular,
\[
 \|u_j\|_{L^2(0,T;L^6)}\le C\|\nabla u_j\|_{L^2((0,T)\times\R^3)} .
\]
Assume:
\begin{enumerate}[label=\textup{(\roman*)},leftmargin=*]
\item The continuity and renormalized continuity equations hold on $\R^3$
after zero extension; the momentum equation holds on $\Omega_j$.
Each equation includes its initial term. On $\Omega_j$, the pressure satisfies
\[
 P_j=a\rho_j^\gamma+\pi_j(\rho_j),\qquad \pi_j\ge0,\qquad
 \sup_{0\le z\le L}\pi_j(z)\longrightarrow0\quad (L<\infty).
\]
\item For each $T>0$ and all $j$ with $T_j>T$,
\begin{equation}\label{cp:uniform_bounds}
 \esssup_{0<t<T}\int_{\R^3}\bigl(\rho_j+\rho_j^\gamma+\rho_j|u_j|^2+P_j\bigr)(t)\dd x
 +\int_0^T\!\!\int_{\R^3}|\nabla u_j|^2\dd x\dd t\le C_T .
\end{equation}
\item There are internal energies $\Theta_j(z)\ge az^\gamma/(\gamma-1)$ such that
\begin{equation}\label{cp:approx_energy}
 \begin{aligned}
 E_j(t)+\int_0^t\mathcal D(u_j(s))\dd s&\le E_{0j}
 &&\text{for a.e. }t,\\
 E_j(t)&=\int_{\Omega_j}\left(\tfrac12 \rho_j|u_j|^2+\Theta_j(\rho_j)\right)(t)\dd x .
 \end{aligned}
\end{equation}
For initial data as in \eqref{eq:initial-data},
\[
 \rho_{0j}\to\rho_0\ \text{in }L^1\cap L^\gamma,\qquad
 m_{0j}\to m_0\ \text{in }L^1,\qquad E_{0j}\to E_0 .
\]
\item For each $I\Subset(0,\infty)$ and $\xi\in C_c^\infty(\R^3)$,
there exists $\sigma=\sigma(I,\xi)>0$ such that, for all sufficiently large $j$,
\begin{equation}\label{cp:momentum_gain}
 \|\xi m_j\|_{L^2(I;H^{-1+\sigma}(\R^3))}\le C_{I,\xi}.
\end{equation}
\end{enumerate}
Then a subsequence converges to a bounded-energy renormalized weak solution
with data $(\rho_0,m_0)$. On each compact interior space--time cylinder,
\[
 \begin{alignedat}{2}
 \rho_j&\to\rho &&\quad\text{strongly in }L^\gamma,\\
 P_j&\to a\rho^\gamma &&\quad\text{strongly in }L^1,\\
 (\rho_j u_j\otimes u_j,\rho_j|u_j|^2)&\rightharpoonup(\rho u\otimes u,\rho|u|^2)
 &&\quad\text{weakly in }L^1 .
 \end{alignedat}
\]
The solution conserves total mass and has the following time representatives:
\begin{equation}\label{eq:weak-traces}
 \rho\in C_{\rm weak}([0,T];L^\gamma),\qquad
 m\in C_{\rm weak}([0,T];L^{2\gamma/(\gamma+1)})
 \qquad (T<\infty).
\end{equation}
Moreover,
\begin{equation}\label{cp:limit_energy}
 E(t)+\int_0^t\mathcal D(u(s))\dd s\le E_0\qquad\text{for a.e. }t>0 .
\end{equation}
Along the full-measure set in \eqref{cp:limit_energy},
\begin{equation}\label{eq:strong-initial-density}
 \rho(t)\longrightarrow\rho_0\quad\text{in }L^1\cap L^\gamma
 \qquad\text{as }t\downarrow0 .
\end{equation}
The initial term in \eqref{eq:renormalized} is given by $b(\rho_0)$.
\end{theorem}

\begin{remark}
The positive margin $\sigma$ has two roles. Combined with the time
derivative bound from the momentum equation, it yields compactness
of the localized momentum in $L^2_tH^{-1}_x$ and permits identification
of the convective flux. It also makes the fractional commutator
pairing strong enough to obtain a positive density moment of the
pressure. Both arguments allow an arbitrarily small positive margin;
neither requires strong $H^{-1}$ convergence of the untruncated density.
The theorem is therefore a conditional compactness statement for
every $\gamma>1$; the range in the existence theorem is determined by
the uniform estimate proved in Section~\ref{ao:section}.

\end{remark}

The main result of this paper states as follows:

\begin{theorem}[Global existence]\label{thm:nine-eighths}
Let $\gamma>9/8$ and assume \eqref{eq:initial-data}. There exists a global
finite-energy renormalized weak solution of \eqref{eq:NS}. It conserves total mass and satisfies
\eqref{eq:weak-traces} and \eqref{eq:strong-initial-density}.
For each $I\Subset(0,\infty)$ and every $\alpha$ satisfying
\begin{equation}\label{intro:signed-gain-range}
 0<2\alpha<\min\left\{\frac16,\frac{8\gamma-9}{4\gamma+1}\right\},
\end{equation}
we have
\begin{equation}\label{intro:signed-gain}
 \|I_{1-\alpha}m\|_{L^2(I\times\R^3)}\le C_{I,\alpha} .
\end{equation}
For $M,A\ge1$,
\begin{equation}\label{intro:signed-tails}
 \begin{aligned}
 \iint_{\{t\in I,\ |u|>M\}}\rho|u|^2\dd x\dd t
    &\le C_{I,\alpha} M^{-\frac{2\alpha}{1+\alpha}},\\
 \iint_{\{t\in I,\ \rho>A\}}\rho|u|^2\dd x\dd t
    &\le C_{I,\alpha} A^{-\frac{(\gamma-1)\alpha}{1+2\alpha}} .
 \end{aligned}
\end{equation}
Estimates \eqref{intro:signed-gain} and \eqref{intro:signed-tails} are
uniform for the artificial-pressure approximations in the proof.
Here $C_{I,\alpha}$ may also depend on the fixed physical parameters
and on the common mass, energy and dissipation bounds used in the
construction. It is independent of the levels $M,A$ and of the
approximation index. The artificial-pressure exponent is fixed
throughout the construction; no uniformity in that exponent is asserted.
\end{theorem}
\begin{remark}
At $\gamma=9/8$, \eqref{intro:signed-gain-range} gives no positive $\alpha$.

Local compactness and the whole-space energy passage require
different inputs. Theorem~\ref{cp:compactness} uses localized momentum
bounds. To pass the differential energy inequality on $\R^3$,
Section~\ref{r4:global-energy} combines the global momentum bound
with a pressure moment, a high-density kinetic-energy tail, and
tightness of the mass. These estimates prevent energy from escaping
to spatial infinity and yield strong convergence of the total energy
on interior time intervals.

\end{remark}

\paragraph{Organization of the paper.}
Section~\ref{sec:compactness} proves Theorem~\ref{cp:compactness}.
Section~\ref{ao:section} proves the momentum and kinetic-energy
estimates \eqref{intro:signed-gain}--\eqref{intro:signed-tails}.
Section~\ref{r4:global-energy} gives the global pressure bound and
energy convergence. Section~\ref{sec:construction} constructs the
approximations and passes to the limit.
At fixed artificial pressure, that construction uses the large
pressure exponent to obtain exact whole-space solutions before the
projection estimate is applied. Appendix~\ref{app:compact-tools}
justifies the pressure tests, Appendix~\ref{app:middle-tools} proves
the projected determinant estimate, and Appendix~\ref{app:momentum-cover}
gives the interval count used in the momentum cover.

\section{Compactness from a fractional momentum estimate}
\label{sec:compactness}

We prove Theorem~\ref{cp:compactness}. We first identify the limits of
$\rho_j u_j$ and $(\rho_j u_j)\otimes u_j$. The effective viscous flux
then identifies the pressure; see \cite{Feireisl2001}.

In this section, constants indexed by $T$ or an interior cylinder may
also depend on the fixed coefficients $a,\gamma,\mu,\lambda$ and the
common bound in \eqref{cp:uniform_bounds} on a containing time interval
$(0,T)$. Constants indexed by an interior
cylinder may additionally depend on the corresponding localized momentum
bound \eqref{cp:momentum_gain} on a larger cylinder. Further dependence
on cutoffs, exponents or truncation levels is indicated by subscripts
or specified locally. These constants are independent of $j$.

\subsection{The weak limit of the momentum}

We truncate the density before taking its product with the weakly
convergent velocity. Choose $\chi_1\in C^\infty([0,\infty))$ nonincreasing, with
$\chi_1=1$ on $[0,1]$ and $\chi_1=0$ on $[3,\infty)$. Set
\begin{equation}\label{cp:density_caps}
 \vartheta_k(z)=\int_0^z\chi_1(s/k)\dd s,\qquad k\in\mathbb N .
\end{equation}
Then $\vartheta_k(z)=z$ for $z\le k$ and $\vartheta_k(z)\uparrow z$ as $k\to\infty$. Moreover,
\[
 0\le \vartheta_k'\le1,\qquad 0\le \vartheta_k(z)-z\vartheta_k'(z)\le \vartheta_k(z)\le z.
\]

\begin{lemma}[Weak limit of $\rho_j u_j$]
\label{cp:identify_momentum}
Let $\gamma>1$. Assume the bounds \eqref{cp:uniform_bounds} for
$\rho_j,\rho_j|u_j|^2,u_j$ and the renormalized continuity equation. Suppose
\[
 \rho_j\rightharpoonup\rho,\quad m_j\rightharpoonup m
       \quad\text{in distributions},\qquad
 u_j\rightharpoonup u\quad\text{in }L^2_{\rm loc}(H^1_{\rm loc}).
\]
Then $\rho u\in L^1_{\rm loc}$ and $m=\rho u$ almost everywhere.
\end{lemma}

\begin{proof}
Fix $k$ and $\chi\in C_c^\infty(\R^3)$. For $f_j=\chi \vartheta_k(\rho_j)$,
the renormalized continuity equation gives
\begin{equation}\label{cp:fixed_cap_equation}
\begin{aligned}
 \partial_t f_j+\diver(f_j u_j)&=g_j,\\
 g_j&=\chi\cdot\bigl(\vartheta_k(\rho_j)-\rho_j\vartheta_k'(\rho_j)\bigr)\diver u_j
           +\vartheta_k(\rho_j)u_j\cdot\nabla\chi .
 \end{aligned}
\end{equation}
On $(0,T)$, for fixed $k$ and $\chi$,
\[
 \|f_j\|_\infty+\|f_j u_j\|_2+\|g_j\|_2
       +\|\partial_t f_j\|_{L^2H^{-1}}\le C_{k,\chi,T}.
\]
By compactness on the support of $\chi$ and \cite{Simon1987}, a diagonal
subsequence satisfies
\[
 \begin{aligned}
 \vartheta_k(\rho_j)&\stackrel{*}{\rightharpoonup}\overline{\vartheta_k(\rho)}
                     &&\text{in }L^\infty,\\
 f_j&\to\chi\overline{\vartheta_k(\rho)} &&\text{in }L^2H^{-1}.
 \end{aligned}
\]
The weak limits are nondecreasing in $k$ and bounded above by $\rho$. Moreover,
\[
 0\leq\rho_j-\vartheta_k(\rho_j)\leq\rho_j{\bf1}_{\{\rho_j>k\}},
 \qquad
 \int_{\R^3}\rho_j{\bf1}_{\{\rho_j>k\}}\,\dd x
       \leq k^{1-\gamma}\int_{\R^3}\rho_j^\gamma\,\dd x .
\]
Hence $\overline{\vartheta_k(\rho)}\uparrow\rho$ a.e. and in $L^1_{\rm loc}$.
At fixed $k$, the $L^2H^{-1},L^2H^1$ pairing gives
\[
 \vartheta_k(\rho_j)u_j\longrightarrow \overline{\vartheta_k(\rho)} u\quad\text{in distributions}.
\]
To remove the truncation, use Cauchy--Schwarz:
\begin{equation}\label{cp:momentum_cap_tail}
\begin{split}
 \int_0^T\!\!\int_{\R^3}
 |(\rho_j-\vartheta_k(\rho_j))u_j|\,\dd x\,\dd t
 &\leq
 \left(\int_0^T\!\!\int_{\{\rho_j>k\}}\rho_j\right)^{1/2}
 \left(\int_0^T\!\!\int_{\R^3}\rho_j|u_j|^2\right)^{1/2}\\
 &\leq C_Tk^{(1-\gamma)/2}.
\end{split}
\end{equation}
Since $m\in L^{2\gamma/(\gamma+1)}_{\rm loc}$, duality with compactly
supported tests and \eqref{cp:momentum_cap_tail} give, on each compact cylinder $Q$,
\[
 \|m-\overline{\vartheta_k(\rho)}u\|_{L^1(Q)}\le C_T k^{(1-\gamma)/2}\longrightarrow0 .
\]
Since $\overline{\vartheta_k(\rho)}u\to\rho u$ almost everywhere, $m=\rho u$.
\end{proof}

\begin{lemma}[The product $m_i u_k$]
\label{cp:kinetic_duality}
Let $\rho\geq0$, let $u\in L^2_{\rm loc}(H^1_{\rm loc})$, and suppose
$m=\rho u\in L^1_{\rm loc}$.
Assume that every compact localization of $m$ belongs to $L^2_tH^{-1}_x$.
Then $\rho|u|^2\in L^1_{\rm loc}$. The $H^{-1},H^1$ pairing of $m_i$
with $u_k$ agrees locally with the function $\rho u_i u_k$.
\end{lemma}

\begin{proof}
For $0\le\zeta\in C_c^\infty$, set
\[
 U_M(v)=
 \begin{cases}
 v,& |v|\leq M,\\
 Mv/|v|,& |v|>M.
 \end{cases}
\]
This truncation satisfies
\[
 |U_M(u)|\le |u|,\qquad
 |\nabla U_M(u)-\nabla u|\le 2|\nabla u|{\bf1}_{\{|u|>M\}}.
\]
Thus $U_M(u)\to u$ in $L^2H^1_{\rm loc}$. For fixed $M$, mollify
$U_M(u)$ near $\operatorname{supp}\zeta$. The approximants are bounded by
$M$ and converge in $L^2H^1_{\rm loc}$ and, along a subsequence, a.e.
Duality and dominated convergence against $\zeta m$ yield
\[
 \langle \zeta m,U_M(u)\rangle
   =\int\zeta\rho|u|^2\min\{1,M/|u|\}\le C_{\zeta,m,u} .
\]
Here the constant depends on the local norms of $m$ and $u$ near
$\supp\zeta$, but is independent of $M$.
Monotone convergence gives $\zeta\rho|u|^2\in L^1$. Now
$|\zeta m_iU_M(u)_k|\le\zeta\rho|u|^2$; dominated convergence proves
the componentwise identity.
\end{proof}

\subsection{Testing the momentum equation with a density truncation}

For the test $\zeta\nabla\Delta^{-1}f$ we use
\[
 \Delta^{-1}=-I_2,\qquad
 \diver(\nabla\Delta^{-1}f)=f,\qquad
 R_{ik}=\partial_i\partial_k\Delta^{-1}\quad (1\le i,k\le3).
\]
We sum repeated indices. The multiplier $\xi_i\xi_k/|\xi|^2$ defines
a bounded operator on $L^p$ ($1<p<\infty$) and on $H^s$ ($s\in\R$).

\begin{lemma}[Estimates for $\nabla\Delta^{-1}$]
\label{cp:potential}
If $w\in H^{-1}(\R^3)$ is supported in a fixed compact set $K$, then
\begin{equation}\label{cp:potential_low}
 \|\nabla\Delta^{-1}w\|_{L^2(\R^3)}\leq C_K\|w\|_{H^{-1}(\R^3)}.
\end{equation}
The adjoint $-\Delta^{-1}\diver$ satisfies the corresponding estimate
for compactly supported vector distributions.
For $3<p<\infty$ and
$g\in L^1(\R^3)\cap L^{2p/(p+2)}(\R^3)$,
\begin{equation}\label{cp:potential_fractional}
 \|\nabla\Delta^{-1}g\|_{H^{1-3/p}}
 \leq C_p\bigl(\|g\|_{L^{2p/(p+2)}}+\|g\|_{L^1}\bigr).
\end{equation}
\end{lemma}

\begin{proof}
For $|\xi|\le1$ test $w$ against $e^{-\mathrm ix\cdot\xi}$ times
a fixed smooth cutoff equal to one near $K$. Then
\[
 |\widehat w(\xi)|\le C_K\|w\|_{H^{-1}},\qquad
 \int_{|\xi|\le1}|\xi|^{-2}\,\dd\xi<\infty .
\]
For $|\xi|\ge1$, the multiplier has modulus $|\xi|^{-1}$.
This proves \eqref{cp:potential_low}, also for the adjoint.
For \eqref{cp:potential_fractional} use
$L^{2p/(p+2)}\hookrightarrow H^{-3/p}$ at high frequencies and
$|\widehat g|\le\|g\|_1$ at low frequencies.
\end{proof}

\begin{lemma}[Fractional velocity commutator]\label{cp:commutator}
For $3<p<\infty$, $v\in H^1(\R^3)$, and $f\in L^p(\R^3)$,
\begin{equation}\label{cp:commutator_estimate}
 \|[v,R_{ik}]f\|_{H^{1-3/p}}
 \leq C_p\|v\|_{H^1}\|f\|_{L^p},\qquad
 [v,R_{ik}]f=vR_{ik}f-R_{ik}(vf).
\end{equation}
\end{lemma}

\begin{proof}
Use \cite[Lemma 2.6]{FLS2016} with $r=s=2$. For $3<p<\infty$,
\[
 \frac12+\frac1p-\frac13<\frac12<\frac12+\frac1p .
\]
The lemma therefore gives regularity $1-3/p$ in $W^{1-3/p,2}=H^{1-3/p}$.
\end{proof}

\begin{lemma}[Test functions obtained from the transport equation]\label{cp:rough_test}
Suppose that the momentum equation holds on an interior cylinder and
that $\rho\ge0$ and $m=\rho u$. On a larger time
interval assume
\[
 \begin{gathered}
 m,\rho u\otimes u,P\in L^\infty_tL^1_x,\qquad S\in L^2_{t,x},\qquad u\in L^2_tH^1_{\rm loc,x},\\
 \xi m\in L^2_tH^{-1+\sigma}_x \quad\text{for local cutoffs }\xi,\qquad 0<\sigma<1.
 \end{gathered}
\]
Let $f\in L^\infty$ and $g\in L^2_{t,x}$ have fixed compact spatial
support and satisfy
\begin{equation}\label{cp:rough_cap_equation}
 \partial_tf+\diver(fu)=g.
\end{equation}
For every smooth cutoff $\zeta$ compactly supported in the cylinder,
testing with $\zeta\nabla\Delta^{-1}f$ gives
\begin{equation}\label{cp:combined_test}
\begin{split}
 \int\zeta\cdot(Pf-S:\nabla^2\Delta^{-1}f)
 ={}&
 -\int(\partial_t\zeta)m\cdot\nabla\Delta^{-1}f\\
 &-\int(\rho u_i u_k-S_{ik})\cdot(\partial_k\zeta)\partial_i\Delta^{-1}f
 -\int P\nabla\zeta\cdot\nabla\Delta^{-1}f\\
 &-\int\zeta m_i[u_k,R_{ik}]f
 -\int\zeta m\cdot\nabla\Delta^{-1}g .
\end{split}
\end{equation}
In the commutator, multiply $u$ by a spatial cutoff equal to one near
the supports of $f$ and $\zeta$. The last two integrals denote Sobolev
dualities; the other integrals are Lebesgue integrals.
\end{lemma}

The proof is in Appendix~\ref{app:flux-test}. In \eqref{cp:combined_test}
the time derivative of the potential and the convective term occur
together in $m_i[u_k,R_{ik}]f$.

\subsection{Proof of the compactness theorem}

\begin{proof}[Proof of Theorem~\ref{cp:compactness}]
Work on $(0,T)$ and exhaust its interior by compact cylinders.
Choose one diagonal subsequence for these cylinders and all truncation levels.

\paragraph{1. Momentum and the convective term.}
The energy and density bounds imply
\begin{equation}\label{cp:momentum_basic}
 \begin{aligned}
 \|m_j(t)\|_{L^{2\gamma/(\gamma+1)}}
 &\leq\|\rho_j(t)\|_{L^\gamma}^{1/2}\|\rho_j(t)|u_j(t)|^2\|_{L^1}^{1/2},\\
 \|m_j(t)\|_{L^1}
 &\leq\|\rho_j(t)\|_{L^1}^{1/2}\|\rho_j(t)|u_j(t)|^2\|_{L^1}^{1/2}.
 \end{aligned}
\end{equation}
Fix $\xi\in C_c^\infty(\R^3)$. Its support lies in $\Omega_j$ for all
large $j$. For $N>5/2$, the two equations and \eqref{cp:uniform_bounds} give
\[
 \|\partial_t(\xi\rho_j)\|_{L^2H^{-N}}
    +\|\partial_t(\xi m_j)\|_{L^2H^{-N}}\le C_{T,\xi,N},
 \qquad T<T_j .
\]
Indeed, $H^N\hookrightarrow W^{1,\infty}$, $\rho_j u_j\otimes u_j,P_j\in L^\infty_tL^1_x$,
and $S_j\in L^2$. Compactness on the support of $\xi$,
\eqref{cp:momentum_gain} and \cite{Simon1987} give
\begin{equation}\label{cp:momentum_strong}
 \xi m_j\longrightarrow\xi m
 \quad\text{in }L^2(I;H^{-1+\sigma'}),\qquad
 I\Subset(0,T),\quad 0<\sigma'<\sigma.
\end{equation}
Reduce $\sigma$ so that $0<\sigma'<\sigma<1$.
Lemma~\ref{cp:identify_momentum} identifies $m=\rho u$. Pair
\eqref{cp:momentum_strong} with $u_j\rightharpoonup u$ in
$L^2H^1_{\rm loc}$ and use Lemma~\ref{cp:kinetic_duality}: in distributions,
\[
 \rho_j u_j\otimes u_j\longrightarrow\rho u\otimes u,\qquad
 \rho_j|u_j|^2\longrightarrow\rho|u|^2\in L^1_{\rm loc}.
\]

To obtain weak convergence in $L^1$, fix a nonnegative smooth cutoff
$\zeta$ and a smooth bounded vector field $v$. Expanding the square gives
\begin{equation}\label{cp:relative_kinetic}
 \int\zeta\rho_j|u_j-v|^2
 \longrightarrow\int\zeta\rho|u-v|^2.
\end{equation}
By H\"older's inequality, for every measurable $E$,
\[
 \int_E\zeta \rho_j|u_j|^2\le 2\int\zeta\rho_j|u_j-v|^2
                     +C_{T,\zeta}\|v\|_\infty^2 |E|^{1-1/\gamma}.
\]
Smooth bounded $v$ are dense in $L^2(\zeta\rho\,\dd x\dd t)$.
Choose $v$ so the limit in \eqref{cp:relative_kinetic} is small, then take
$|E|$ small. This gives uniform integrability for all large $j$; each
of the finitely many remaining $\rho_j|u_j|^2$ is integrable.
Since $|\rho_j u_j\otimes u_j|\le \rho_j|u_j|^2$, this gives
\[
 \rho_j u_j\otimes u_j\rightharpoonup\rho u\otimes u,\qquad
 \rho_j|u_j|^2\rightharpoonup\rho|u|^2
 \qquad\text{weakly in }L^1_{\rm loc}.
\]
Testing in time and exhausting $\R^3$ gives
\[
 \rho,\ \rho|u|^2,\ m=\rho u\ \in L^\infty(0,T;L^1(\R^3)).
\]
\paragraph{2. Higher integrability of the pressure.}
Choose cutoffs $0\le\zeta,\chi\le1$ with $\zeta=1$ on the
interior cylinder $Q$ and $\chi=1$ on the spatial support of $\zeta$.
For a smooth nonincreasing cutoff $\chi_0$ equal to one on $[0,1]$
and zero on $[2,\infty)$, set
\[
 \frac3p<\sigma',\qquad 0<\theta<1,\qquad \theta p\le\gamma,
 \qquad b_N(z)=\int_0^z\theta(1+s)^{\theta-1}\chi_0(s/N)\,\dd s.
\]
With these choices fixed, $C_{Q,p,\theta}$ denotes constants that may
depend on $p,\theta$, the supports and smooth norms of the chosen
cutoffs, and the common energy and localized momentum bounds on the
larger cylinder, but not on $j$ or $N$.
For fixed $N$, $b_N$ is bounded with compactly supported derivative.
As $N\to\infty$, $b_N(z)\uparrow(1+z)^\theta-1$, and
\[
 |b_N(z)|+|zb_N'(z)|\leq 2z^\theta
\]
uniformly in $N$.
For $f_j=\chi b_N(\rho_j)$, renormalized continuity reads
\begin{equation}\label{cp:power_cap_equation}
\begin{split}
 \partial_tf_j+\diver(f_ju_j)&=g_j,\\
 g_j&=\chi\cdot(b_N(\rho_j)-\rho_jb_N'(\rho_j))\diver u_j
       +b_N(\rho_j)u_j\cdot\nabla\chi.
\end{split}
\end{equation}
H\"older's inequality and the local Sobolev bound on $u_j$ give
\[
 \sup_{j,N}\bigl(\|f_j\|_{L^\infty_tL^p_x}
             +\|g_j\|_{L^2_tL^{2p/(p+2)}_x}\bigr)<\infty.
\]
For fixed $N$, $g_j\in L^2$ since $b_N-zb_N'$ is bounded.
Lemma~\ref{cp:rough_test} applies. With the same velocity cutoff,
Lemmas~\ref{cp:potential}--\ref{cp:commutator} give, uniformly in $j,N$,
\[
 \|[u_{j,k},R_{ik}]f_j\|_{L^2H^{1-3/p}}
   +\|\nabla\Delta^{-1}g_j\|_{L^2H^{1-3/p}}\le C_{Q,p,\theta} .
\]
These terms pair with $\xi m_j$ since $3/p<\sigma'$. For the stress
and the derivatives of $\zeta$, the compact support of $f_j$ and $p>3$ give
\[
 \|R_{ik}f_j\|_{L^\infty_tL^2_x}
    +\|\nabla\Delta^{-1}f_j\|_{L^\infty_{t,x}}\le C_{Q,p,\theta} .
\]
Thus \eqref{cp:combined_test}, \eqref{cp:uniform_bounds}
and \eqref{cp:momentum_basic} imply
\[
 0\le \int\zeta P_j f_j\le C_{Q,p,\theta} .
\]
Let $N\to\infty$ by monotone convergence. Using the $L^1$ pressure
bound on $\{\rho_j\le 1\}$, we obtain
\begin{equation}\label{cp:pressure_moment}
 \int_QP_j\rho_j^\theta\,\dd x\,\dd t\leq C_{Q,p,\theta} .
\end{equation}
Since $P_j\ge a\rho_j^\gamma$, the family $\{\rho_j^\gamma\}$ is uniformly
integrable on $Q$.
For $L\ge1$,
\[
 \int_Q\pi_j(\rho_j)
 \leq |Q|\sup_{0\leq r\leq L}\pi_j(r)
       +L^{-\theta}\int_QP_j\rho_j^\theta.
\]
Let $j\to\infty$ first and then $L\to\infty$. After extraction,
\[
 \pi_j(\rho_j)\to0\quad\text{in }L^1_{\rm loc},\qquad
 P_j\rightharpoonup\overline{p(\rho)}
                  \quad\text{weakly in }L^1_{\rm loc}.
\]
By exhaustion $\overline{p(\rho)}\in L^\infty(0,T;L^1)$.
The limiting momentum equation has pressure $\overline{p(\rho)}$.

\paragraph{3. Weak continuity of the effective viscous flux.}
Fix a smooth bounded $b$ with compactly supported derivative.
Set $f_j=\chi b(\rho_j)$ and define $g_j$ as in
\eqref{cp:power_cap_equation} with $b$ in place of $b_N$. Then
\[
 f_j\to\bar f=\chi\overline{b(\rho)}\quad\text{in }L^2H^{-1},
 \qquad g_j\rightharpoonup\bar g\quad\text{in }L^2 .
\]
Here $\overline{b(\rho)}$ denotes the weak limit, which has not yet
been identified with $b(\rho)$. Pairing with $u_j\rightharpoonup u$
in local $L^2H^1$ gives
\[
 f_ju_j\to\bar fu \quad\text{in distributions},\qquad
 \partial_t\bar f+\diver(\bar f u)=\bar g.
\]

Since $R_{ik}f_j\to R_{ik}\bar f$ in $L^2H^{-1}$, the same pairing
on the support of the velocity cutoff gives
\[
 u_{j,k}R_{ik}f_j\to u_kR_{ik}\bar f,\qquad
 u_{j,k}f_j\to u_k\bar f\qquad\text{in distributions}.
\]
For fixed $b$, $u_{j,k}f_j$ is bounded in $L^2_{t,x}$: $f_j$ is
bounded with fixed compact spatial support. Its distributional limit is
therefore also its weak $L^2$ limit. Boundedness of $R_{ik}$ on $L^2$
identifies the limit of $R_{ik}(u_{j,k}f_j)$. Lemma~\ref{cp:commutator} gives
\[
 [u_{j,k},R_{ik}]f_j\rightharpoonup[u_k,R_{ik}]\bar f
          \quad\text{in }L^2H^{1-3/p}.
\]
This weak convergence pairs with \eqref{cp:momentum_strong} since
$3/p<\sigma'$. For the source term Lemma~\ref{cp:potential} gives
\[
 -\Delta^{-1}\diver(\zeta m_j)\to-\Delta^{-1}\diver(\zeta m)
       \quad\text{in }L^2_{t,x}.
\]
Pair this limit with $g_j\rightharpoonup\bar g$ in $L^2$.
For the terms with derivatives of $\zeta$, Lemma~\ref{cp:potential}
and interpolation yield
\[
 \begin{aligned}
 \nabla\Delta^{-1}f_j&\to\nabla\Delta^{-1}\bar f
                         &&\text{in }L^2_tL^2_x,\\
 \sup_j\|\nabla\Delta^{-1}f_j\|_{L^\infty_tW^{1,p}(K)}&<\infty
                         &&(3<p<\infty,\ K\text{ a ball}),\\
 \nabla\Delta^{-1}f_j&\to\nabla\Delta^{-1}\bar f
                         &&\text{in }L^2_tL^\infty_{\rm loc,x}.
 \end{aligned}
\]
Thus we may replace the potential by its limit using the time-slice
$L^1$ bounds, then use weak $L^1_{\rm loc}$ convergence of $m_j,\rho_j u_j\otimes u_j,P_j$
and weak $L^2$ convergence of $S_j$.

For the remaining stress term, integration by parts gives
\begin{equation}\label{cp:stress_localization}
 \int\varphi S_j:\nabla^2\Delta^{-1}f_j
 =(2\mu+\lambda)\int\varphi\cdot(\diver u_j)f_j
   +2\mu\int u_{j,i}\cdot
        \bigl(\partial_i\varphi f_j-\partial_k\varphi R_{ik}f_j\bigr).
\end{equation}
The terms with $\nabla\varphi$ converge by $L^2H^{-1},L^2H^1$
duality. Subtract \eqref{cp:combined_test} for the limit fields
from the limit of \eqref{cp:combined_test} for the approximations. On
$\{\chi=1\}$, Lemma~\ref{cp:rough_test} yields
\begin{equation}\label{cp:effective_flux}
 \overline{a\rho^\gamma b(\rho)}
       -\overline{p(\rho)}\,\overline{b(\rho)}
 =(2\mu+\lambda)\left(
       \overline{b(\rho)\diver u}
       -\overline{b(\rho)}\,\diver u\right)
 \quad\hbox{in distributions}.
\end{equation}
The artificial pressure makes no contribution since $b$ is bounded
and $\pi_j(\rho_j)\to0$ in local $L^1$.

\paragraph{4. Oscillations of the density.}
Let $X$ have the Young measure of $\rho_j$ at $(t,x)$; write
$\mathbb E$ for integration against this measure. By uniform
integrability (Step~2), for bounded continuous $b$,
\[
 \rho=\mathbb EX,\qquad
 \overline{p(\rho)}=a\mathbb EX^\gamma,\qquad
 \overline{p(\rho)b(\rho)}=a\mathbb E[X^\gamma b(X)] .
\]

For a bounded nondecreasing $1$-Lipschitz function $b$, put
\[
 \Gamma_b=a\bigl(\mathbb E[X^\gamma b(X)]-(\mathbb E X^\gamma)\mathbb E b(X)\bigr).
\]
If $Y$ is an independent copy of $X$, then

\begin{equation}\label{cp:covariance}
\begin{split}
 \Gamma_b
 &=\frac a2\mathbb E[(X^\gamma-Y^\gamma)(b(X)-b(Y))]\\
 &\geq\frac a2\mathbb E|b(X)-b(Y)|^{\gamma+1}
 \geq\frac a2\mathbb E|b(X)-\mathbb Eb(X)|^{\gamma+1}\geq0 .
\end{split}
\end{equation}
Here we used $|x^\gamma-y^\gamma|\ge|x-y|^\gamma$ for $x,y\ge0$
and then Jensen's inequality with $X$ fixed.

Take $b=\vartheta_k$, so $\overline{\vartheta_k(\rho)}=\mathbb E \vartheta_k(X)$. The right-hand side
of \eqref{cp:effective_flux} is the distributional limit of
\[
 (2\mu+\lambda)(\diver u_j-\diver u)\cdot
       (\vartheta_k(\rho_j)-\overline{\vartheta_k(\rho)}).
\]
This follows by expanding the product and using $\overline{\vartheta_k(\rho)}\in L^\infty$
and $\diver u\in L^2$. In this step, $C_T$ is independent of the
truncation level $k$ and the spatial exhaustion cutoff $\zeta$;
it depends only on the fixed coefficients and the common mass and
dissipation bounds on $(0,T)$. For $0\le\zeta\le1$ compactly supported,
Cauchy--Schwarz and the global dissipation bound give
\begin{equation}\label{cp:oscillation_first}
 \int\zeta\,\mathbb E|\vartheta_k(X)-\overline{\vartheta_k(\rho)}|^{\gamma+1}
 \leq C_T
      \left(\int\zeta\,\mathbb E|\vartheta_k(X)-\overline{\vartheta_k(\rho)}|^2\right)^{1/2}.
\end{equation}
For fixed $k$, all centered powers here are bounded.
H\"older's inequality on the product measure gives
\begin{equation}\label{cp:oscillation_interpolation}
\int\zeta\,\mathbb E|\vartheta_k(X)-\overline{\vartheta_k(\rho)}|^2
 \le
 \left(\int\zeta\,\mathbb E|\vartheta_k(X)-\overline{\vartheta_k(\rho)}|\right)^{1-1/\gamma}
 \left(\int\zeta\,\mathbb E|\vartheta_k(X)-\overline{\vartheta_k(\rho)}|^{\gamma+1}\right)^{1/\gamma}.
\end{equation}
The first integral on the right is at most
$2\int_0^T\!\int\rho$. Combining the last two inequalities and using
$1/(2\gamma)<1$ bounds the left-hand side of
\eqref{cp:oscillation_first} uniformly in $k$ and $\zeta$.
Exhaust $(0,T)\times\R^3$ to obtain

\begin{equation}\label{cp:centered_moments}
 \sup_k\int_0^T\!\!\int_{\R^3}
 \left(\mathbb E|\vartheta_k(X)-\overline{\vartheta_k(\rho)}|^{\gamma+1}
       +\mathbb E|\vartheta_k(X)-\overline{\vartheta_k(\rho)}|^2\right)\,\dd x\,\dd t
 \leq C_T.
\end{equation}
By Fatou's lemma and $\overline{\vartheta_k(\rho)}\uparrow\rho$,
\begin{equation}\label{cp:full_centered_moments}
 \int_0^T\!\!\int_{\R^3}
 \left(\mathbb E|X-\rho|^{\gamma+1}
       +\mathbb E|X-\rho|^2\right)\,\dd x\,\dd t\leq C_T .
\end{equation}
In particular, $\mathbb E|X-\rho|\in L^2$ by Jensen's inequality;
no $L^2$ estimate on $\rho$ is needed.

\paragraph{5. The renormalized continuity equation and its initial value.}
On $-1<t<0$ set $(\rho_j,u_j)=(\rho_{0j},0)$.
The initial terms extend both continuity equations across $t=0$.
Strong convergence of $\rho_{0j}$ gives the limits $\rho_0$ and
$\vartheta_k(\rho_0)$ on this interval. Then
\begin{equation}\label{cp:limit_truncation_equation}
 \begin{aligned}
 \partial_t \overline{\vartheta_k(\rho)}
 +\diver(\overline{\vartheta_k(\rho)}u)&=G_k,\\
 G_k&=\overline{(\vartheta_k(\rho)-\rho \vartheta_k'(\rho))\diver u}.
 \end{aligned}
\end{equation}
Here $G_k\in L^2$ on each finite time interval and vanishes for $t<0$.
To renormalize this equation, use a fixed standard mollifier profile
$J_1$ and set $J_\varepsilon(x)=\varepsilon^{-3}J_1(x/\varepsilon)$:
\[
 \begin{aligned}
 \overline{\vartheta_k(\rho)}^\varepsilon
 &=J_\varepsilon*\overline{\vartheta_k(\rho)},\\
 r_\varepsilon
 &=\diver(\overline{\vartheta_k(\rho)}^\varepsilon u)
      -J_\varepsilon*\diver(\overline{\vartheta_k(\rho)}u).
 \end{aligned}
\]
On a fixed compact cylinder,
\[
 r_\varepsilon(x)
 =(\diver u)(x)\overline{\vartheta_k(\rho)}^\varepsilon(x)
  +\int\nabla J_\varepsilon(z)\cdot
       [u(x)-u(x-z)]\overline{\vartheta_k(\rho)}(x-z)\,\dd z .
\]
On nested compact cylinders $Q\Subset Q'$, for all sufficiently small
$\varepsilon$ depending on the gap between the cylinders, the translation
estimate gives
\[
 \|r_\varepsilon\|_{L^2(Q)}
       \le C_{J_1}\|\overline{\vartheta_k(\rho)}\|_\infty\|\nabla u\|_{L^2(Q')}.
\]
For smooth $u$ the two terms cancel as $\varepsilon\downarrow0$ by
Taylor's formula, translation continuity of $\overline{\vartheta_k(\rho)}$ in $L^2_{\rm loc}$, and
\[
 \int \partial_iJ_\varepsilon(z) z_\ell\,\dd z=-\delta_{i\ell}.
\]
Approximating $u$ in $L^2H^1_{\rm loc}$ therefore gives
$r_\varepsilon\to0$ in local $L^2$; see \cite{DiPernaLions1989}.
For fixed $\varepsilon$, $\partial_t \overline{\vartheta_k(\rho)}^\varepsilon\in L^2_{\rm loc}$.
Use the chain rule for smooth $h$ and let $\varepsilon\downarrow0$:
\begin{equation}\label{cp:renormalized_tk}
 \begin{split}
 \partial_t h(\overline{\vartheta_k(\rho)})
 +\diver(h(\overline{\vartheta_k(\rho)})u)
 &\\
 +\bigl(\overline{\vartheta_k(\rho)}h'(\overline{\vartheta_k(\rho)})
        -h(\overline{\vartheta_k(\rho)})\bigr)\diver u
 &=h'(\overline{\vartheta_k(\rho)})G_k .
 \end{split}
\end{equation}
locally, including across the extended initial time.

Suppose $h'=0$ on $[M,\infty)$ and $k>2M$. The coefficient of
$\diver u_j$ in $G_k$ vanishes for $\rho_j\le k$. On
$\{\overline{\vartheta_k(\rho)}\le M,\ \rho_j>k\}$,
\[
 \vartheta_k(\rho_j)-\overline{\vartheta_k(\rho)}\ge k/2,\qquad
 \vartheta_k(\rho_j)\le2|\vartheta_k(\rho_j)-\overline{\vartheta_k(\rho)}|.
\]
Hence
\[
 {\bf1}_{\{\overline{\vartheta_k(\rho)}\leq M\}}
 |\vartheta_k(\rho_j)-\rho_j\vartheta_k'(\rho_j)|^2
 \leq4(k/2)^{1-\gamma}|\vartheta_k(\rho_j)-\overline{\vartheta_k(\rho)}|^{\gamma+1}.
\]
At fixed $k$, simple-function approximation of $\overline{\vartheta_k(\rho)}$ identifies the
weak limits of these bounded centered powers. Cauchy--Schwarz and
\eqref{cp:centered_moments} now give locally
\begin{equation}\label{cp:renormalization_source}
 \|h'(\overline{\vartheta_k(\rho)})G_k\|_{L^1}
 \leq C_{h,T}k^{(1-\gamma)/2}\longrightarrow0.
\end{equation}
As $\overline{\vartheta_k(\rho)}\uparrow\rho$, both
\[
 h(\overline{\vartheta_k(\rho)})\quad\text{and}\quad
 \overline{\vartheta_k(\rho)}h'(\overline{\vartheta_k(\rho)})
       -h(\overline{\vartheta_k(\rho)})
\]
remain bounded. In \eqref{cp:renormalized_tk}, dominated convergence
and \eqref{cp:renormalization_source} give \eqref{eq:renormalized}
with initial value $h(\rho_0)$. Approximation in $C^1$, with all
derivatives supported in a fixed interval, gives the same identity
for every admissible $h\in C^1$.

\paragraph{6. Strong convergence of the density.}
On $\R^3$, finite mass does not ensure integrability of the negative
part of $\rho\log\rho$. Choose a smooth nondecreasing
$\chi_2:[0,\infty)\to[0,1]$, zero on $[0,1]$ and one on $[2,\infty)$,
and define
\[
 \widetilde{\vartheta}_k(z)=\int_0^z\vartheta_k'(s)\chi_2(ks)\,\dd s,\qquad
 \ell_k''(z)=\frac{\widetilde{\vartheta}_k'(z)}{z},\qquad
 \ell_k(0)=\ell_k'(0)=0.
\]
Since $\widetilde{\vartheta}_k'=0$ near zero,
\begin{equation}\label{cp:entropy_properties}
\begin{gathered}
 0\leq \widetilde{\vartheta}_k'\leq1,\qquad 0\leq \widetilde{\vartheta}_k(z)\leq z,\qquad \widetilde{\vartheta}_k(z)\to z,\\
 \ell_k\text{ is convex},\qquad
 0\leq \ell_k(z)\leq C_kz,\qquad
 z\ell_k'(z)-\ell_k(z)=\widetilde{\vartheta}_k(z).
\end{gathered}
\end{equation}
Write $\ell_k(z)=c_k z+[\ell_k(z)-c_k z]$, where $c_k=\ell_k'(z)$ for
$z\ge3k$. The bracket has compactly supported derivative: use
renormalized continuity for this part and continuity for $c_k\rho$.
For each fixed $k$, the bracket is bounded and compact in local
$L^2H^{-1}$ along the sequence, as in Step~3. The linear flux converges
by $m=\rho u$. These limits also hold in the entropy equation tested
across $t=0$: $|\ell_k(\rho_j)u_j|\le C_k|m_j|$ and
$|\widetilde{\vartheta}_k(\rho_j)\diver u_j|\le C_k|\diver u_j|$.
On $(0,\varepsilon)\times K$ their integrals are
$O_{k,T}(\varepsilon)$ and $O_{k,K,T}(\sqrt\varepsilon)$, respectively.
Moreover $\ell_k$ is Lipschitz, so
$\ell_k(\rho_{0j})\to \ell_k(\rho_0)$ in $L^1$. The limiting entropy
equation and that of $\rho$ thus have the same initial term.
The bound $\ell_k(z)\le C_kz$ and convexity give
\[
 d_k=\mathbb E\ell_k(X)-\ell_k(\rho)\ge0 .
\]
Subtract the entropy equation for $\rho$ from the limit equation,
and apply \eqref{cp:effective_flux} with $b=\widetilde{\vartheta}_k$:
\begin{equation}\label{cp:entropy_defect}
 \partial_td_k+\diver(d_ku)+(2\mu+\lambda)^{-1}\Gamma_{\widetilde{\vartheta}_k}
 =(\diver u)\cdot\bigl(\widetilde{\vartheta}_k(\rho)-\mathbb E\widetilde{\vartheta}_k(X)\bigr).
\end{equation}
Here \eqref{cp:effective_flux} gives
\[
 \overline{\widetilde{\vartheta}_k(\rho)\diver u}-\widetilde{\vartheta}_k(\rho)\diver u
 =\frac{\Gamma_{\widetilde{\vartheta}_k}}{2\mu+\lambda}
    +(\mathbb E\widetilde{\vartheta}_k(X)-\widetilde{\vartheta}_k(\rho))\diver u .
\]
For $t<0$, $d_k=0$. On $(0,T)$,
\[
 0\le d_k\le C_k\rho,\qquad |d_k u|\le C_k|m|,\qquad
 0\le \Gamma_{\widetilde{\vartheta}_k}\le2a\|\widetilde{\vartheta}_k\|_\infty\mathbb EX^\gamma .
\]
For the right-hand side use the Lipschitz bound
\[
 |\widetilde{\vartheta}_k(\rho)-\mathbb E\widetilde{\vartheta}_k(X)|
 \leq\mathbb E|X-\rho|\in L^2
\]
and $\diver u\in L^2$. Integrate \eqref{cp:entropy_defect} against
$\chi_R\uparrow1$, $|\nabla\chi_R|\le C/R$, from a negative time
to $t$. The preceding bounds give $\int\chi_R d_k\in W^{1,1}(-1,T)$
and a transport error at most $C_k\|m\|_{L^1_{t,x}}/R$.
Let $R\to\infty$. For a.e. $t$,
\begin{equation}\label{cp:entropy_integrated}
\begin{split}
 \int_{\R^3}d_k(t)\,\dd x
 +(2\mu+\lambda)^{-1}\int_0^t\!\!\int_{\R^3}\Gamma_{\widetilde{\vartheta}_k}\,\dd x\,\dd s
 \leq{}&
 \|\diver u\|_{L^2((0,T)\times\R^3)}\\
 &{}\times
 \|\widetilde{\vartheta}_k(\rho)-\mathbb E\widetilde{\vartheta}_k(X)\|_{L^2((0,T)\times\R^3)}.
\end{split}
\end{equation}
Choose $t$ in the common full-measure set for all $k$. Since $\widetilde{\vartheta}_k(z)\to z$,
\[
 \widetilde{\vartheta}_k(\rho)-\mathbb E\widetilde{\vartheta}_k(X)\longrightarrow0
 \quad\text{in }L^2((0,T)\times\R^3)
\]
by dominated convergence and \eqref{cp:full_centered_moments}.
Drop $d_k\ge0$ and use \eqref{cp:covariance} and Fatou's lemma:
\[
 \int_0^t\!\!\int_{\R^3}\mathbb E|X-\rho|^{\gamma+1}\,\dd x\,\dd s=0.
\]
Let $t\uparrow T$ in that set. The Young measure is $\delta_{\rho(t,x)}$
a.e., hence $\rho_j\to\rho$ locally in measure. By
\eqref{cp:pressure_moment}, Vitali's theorem and Step~2,
\[
 \rho_j\to\rho\ \text{in }L^\gamma_{\rm loc},\qquad
 a\rho_j^\gamma\to a\rho^\gamma,\quad P_j\to a\rho^\gamma
           \ \text{in }L^1_{\rm loc}.
\]

\paragraph{7. Initial data, total mass and the energy inequality.}
For a test meeting $t=0$, split the time integrals at $\varepsilon$.
On its compact spatial support $K$, uniformly in $j$,
\[
 \int_0^\varepsilon\!\int_K(\rho_j+|m_j|+|\rho_j u_j\otimes u_j|+P_j)\le C_T\varepsilon,
 \qquad
 \int_0^\varepsilon\!\int_K |S_j|\le C_{T,K}\varepsilon^{1/2}.
\]
The same bounds hold for the limit. Let $j\to\infty$ and then
$\varepsilon\downarrow0$. The prepared data give \eqref{eq:weak-mass}
and \eqref{eq:weak-momentum}; Step~5 gives \eqref{eq:renormalized}.
For spatial cutoffs $\chi_R\uparrow1$ as above,
\[
 \int_0^T\!\int |m\cdot\nabla\chi_R|
 \le C R^{-1}\|m\|_{L^1_{t,x}}\longrightarrow0.
\]
Thus mass is conserved a.e. in time; Proposition~\ref{cp:time_traces}
gives the identity for every time.

To pass the energy inequality, take $w\geq0$ smooth and compactly
supported in $(0,T)$.
Fubini in \eqref{cp:approx_energy} gives
\begin{equation}\label{cp:weighted_energy}
 \int_0^T w(t)E_j(t)\,\dd t
 +\int_0^T\left(\int_s^T w(t)\,\dd t\right)\mathcal D(u_j(s))\,\dd s
 \leq E_{0j}\int_0^T w(t)\,\dd t.
\end{equation}
By local convergence of the pressure and kinetic energy, exhaustion
in space and $\Theta_j(z)\ge az^\gamma/(\gamma-1)$,
\[
 \int_0^T w E\le\liminf_j\int_0^T wE_j .
\]
The dissipation is weakly lower semicontinuous with the nonnegative
weight $\int_s^T w(t)\,\dd t$
since $\mu+\lambda\ge\mu/3$. Let $j\to\infty$ in \eqref{cp:weighted_energy}.
As $E_{0j}\to E_0$ and $w\ge0$ is arbitrary, \eqref{cp:limit_energy}
holds a.e. in $t$. A diagonal exhaustion gives the global subsequence.
The initial trace follows from the next proposition.
\end{proof}

\subsection{Weak continuity in time and the initial density}

\begin{proposition}\label{cp:time_traces}
The solution above has the weakly continuous density and momentum
representatives stated in Theorem~\ref{cp:compactness}.
Total mass is conserved at every time.
As $t\downarrow0$ through the full-measure set on which the energy
inequality holds,
\begin{equation}\label{cp:strong_initial_density}
 \rho(t)\longrightarrow\rho_0
 \quad\hbox{strongly in }L^1(\R^3)\cap L^\gamma(\R^3).
\end{equation}
\end{proposition}

\begin{proof}
For $\varphi,\Phi$ smooth and compactly supported in space,
\[
 \begin{aligned}
 \frac{\mathrm d}{\mathrm dt}\int\rho\varphi
     &=\int m\cdot\nabla\varphi &&\in L^\infty(0,T),\\
 \frac{\mathrm d}{\mathrm dt}\int m\cdot\Phi
     &=\int\big[(\rho u\otimes u-S(u)):\nabla\Phi
                    +p(\rho)\diver\Phi\big] &&\in L^1(0,T).
 \end{aligned}
\]
The weak formulations prescribe the values at $t=0$.
Density in the dual spaces and the energy bounds yield
\[
 \rho\in C_{\rm weak}([0,T];L^\gamma),\qquad
 m\in C_{\rm weak}([0,T];L^{2\gamma/(\gamma+1)}).
\]
Positivity and the norm bounds extend to every time by weak lower
semicontinuity. Testing the continuity equation with $\chi_R$ gives
\begin{equation}\label{cp:mass_tightness}
 \left|\int\chi_R\rho(t)-\int\chi_R\rho_0\right|
 \leq\frac{CT}{R}\|m\|_{L^\infty(0,T;L^1)}.
\end{equation}
Monotone exhaustion proves $\int\rho(t)=\int\rho_0$ for every $t$.
It also yields the uniform tail estimate
\[
 \int(1-\chi_R)\rho(t)
 \leq\int(1-\chi_R)\rho_0
       +\frac{CT}{R}\|m\|_{L^\infty(0,T;L^1)}.
\]

The kinetic energy is convex in $(\rho,m)$ and lower semicontinuous for
strong convergence by Fatou's lemma, with the vacuum convention. Hence it
is weakly lower semicontinuous on $L^\gamma\times L^{2\gamma/(\gamma+1)}$.
Along times satisfying \eqref{cp:limit_energy},
\[
\begin{aligned}
 \frac a{\gamma-1}\limsup_{t\downarrow0}\|\rho(t)\|_\gamma^\gamma
 &\le E_0-\liminf_{t\downarrow0}
          \int_{\R^3}\frac{|m(t)|^2}{2\rho(t)}\,\dd x \\
 &\le E_0-\int_{\R^3}\frac{|m_0|^2}{2\rho_0}\,\dd x
   =\frac a{\gamma-1}\|\rho_0\|_\gamma^\gamma.
\end{aligned}
\]
Weak continuity gives the reverse lower bound. Uniform convexity of
$L^\gamma$ yields strong $L^\gamma$ convergence. This gives $L^1$
convergence on balls; the uniform tail estimate gives convergence in
$L^1(\R^3)$.
\end{proof}

\begin{remark}\label{cp:scope}
The assumption $\gamma>1$ gives the integrability exponent
$2\gamma/(\gamma+1)>1$, the tail in \eqref{cp:momentum_cap_tail}
and the decay in \eqref{cp:renormalization_source}.
The additional gain $\sigma>0$ permits $3/p<\sigma'<\sigma$ in the
pressure test. No compactness assumption on $\rho_j$ in $H^{-1}$
is imposed.
\end{remark}

\section{The momentum estimate for \texorpdfstring{$\gamma>9/8$}{gamma > 9/8}}
\label{ao:section}
We prove the fractional estimate required by Theorem~\ref{cp:compactness}.
We cover the large values of a spatial average of $m$ using the momentum
equation. The selected time intervals remain active even at points
where the signed average of the momentum cancels. The full averaged
kinetic tensor retains the corresponding velocity covariance: equal
density-weighted contributions with opposite velocities have zero
average momentum but nonzero kinetic energy. The interval count uses the flux only
during the selected times. The covariance gives a determinant lower bound
on the entire selected cover (Lemmas~\ref{ao:active-product} and
\ref{ao:minor}). Section~\ref{ao:feedback} combines the two estimates.

Let \(J\) be a bounded interval and fix \(I=(t_-,t_+)\Subset J\).
Choose \(0\le\zeta\in C_c^\infty(J)\) with \(\zeta=1\) on \(\overline I\).
For $\rho\ge0$, set $m=\rho u$.
With \(P\) the total pressure and \(S=S(u)\) as in \eqref{eq:stress},
assume on \(J\times\R^3\) that
\begin{equation}\label{ao:balance}
 \begin{gathered}
 \partial_t\rho+\operatorname{div}m=0,\\
 \partial_tm+\operatorname{div}(\rho u\otimes u+P\mathrm{Id}_3)=\operatorname{div}S,
 \qquad P\ge a\rho^\gamma,
 \end{gathered}
\end{equation}
and
\begin{equation}\label{ao:natural}
 \begin{aligned}
 &\operatorname*{ess\,sup}_{t\in J}
 \int_{\R^3}(\rho+\rho^\gamma+\rho|u|^2+P)\,\mathrm dx
 \\
 &\qquad+\|\nabla u\|_{L^2(J\times\R^3)}^2
 +\|S\|_{L^2(J\times\R^3)}^2\le C_0.
 \end{aligned}
\end{equation}
Here $C_0$ is a fixed common bound. In this section, $C_*$ may depend
on $a,\gamma,C_0,I,J$, the selected time cutoff and the fixed reference
kernels and cutoffs. It is independent of $r,h,\Lambda,M,A$ and of
approximation indices. Any additional dependence on $\alpha$ or on a
higher-integrability exponent is displayed separately.
We use the \(L^6\) representative of \(u(t)\in\dot H^1(\R^3)\).
The balance \eqref{ao:balance} is required on the whole space:
zero extension from a ball can introduce a boundary term.

The mass and energy bounds also give
\[
 \begin{gathered}
 m\in L^\infty\!\left(J;L^1(\R^3)
                  \cap L^{\frac{2\gamma}{\gamma+1}}(\R^3)\right),\\
 \rho u\otimes u+P\mathrm{Id}_3\in L^\infty(J;L^1(\R^3)).
 \end{gathered}
\]
\begin{proposition}[Fractional estimate for the momentum]\label{ao:gain}
Suppose \(\gamma>9/8\) and \eqref{ao:balance}--\eqref{ao:natural} hold.
For each \(I\Subset J\),
\[
 \|I_{1-\alpha}(\rho u)\|_{L^2(I\times\R^3)}
 \le C_{I,\alpha},
 \qquad
 0<2\alpha<
 \min\left\{\frac16,\frac{8\gamma-9}{4\gamma+1}\right\}.
\]
For the same \(\alpha\),
\[
 \begin{aligned}
 \iint_{\substack{t\in I\\ |u|>M}}\rho|u|^2
 &\le C_{I,\alpha} M^{-\frac{2\alpha}{1+\alpha}},\qquad M\ge1,\\
 \iint_{\substack{t\in I\\ \rho>A}}\rho|u|^2
 &\le C_{I,\alpha} A^{-\frac{(\gamma-1)\alpha}{1+2\alpha}},\qquad A\ge1 .
 \end{aligned}
\]
With $J$ fixed, the constants $C_{I,\alpha}$ may also depend on
$J,a,\gamma,C_0$ and the fixed reference kernels. The time cutoff is chosen in terms of $I\Subset J$.
They are independent of $M,A$ and are uniform for approximation families
satisfying the same bound $C_0$.
\end{proposition}

\subsection{A projection estimate for the regularized equations}
\label{ao:smoothing}

For $\omega\in\mathbb S^3$, let $\Pi_\omega$ denote orthogonal projection
onto $\omega^\perp$, expressed in an orthonormal basis. Integration along
the lines parallel to $\omega$ gives
\[
 X_\omega f(y)=\int_{\mathbb R}f(y+s\omega)\dd s,\qquad y\in\omega^\perp .
\]
On $\omega^\perp\simeq\R^3$, choose a smooth radial kernel $J_1$
with $J_1\ge0$, $\supp J_1\subset B_1$, and $\int J_1=1$, and set
$J_h(y)=h^{-3}J_1(y/h)$. Define
\[
 \mathcal A_\omega=\Pi_\omega (X_\omega\mathcal A)\Pi_\omega^{\top},
 \qquad
 \mathcal G_h(\mathcal A)=\int_{\mathbb S^3}\int_{\omega^\perp}
        \det(J_h*\mathcal A_\omega)^{1/2}\dd y\dd\omega .
\]

\begin{proposition}[Upper bound for the projected determinant]\label{rt:optimized-CI}
Let $\mathcal A$ and $\mathcal E$ be matrix fields on $\R^4$ such that
\[
 \mathcal A=\mathcal A^{\top}\ge0,\qquad \mathcal A\in L^1,\qquad
 \mathcal E\in L^2,\qquad \Div(\mathcal A-\mathcal E)=\mathcal F\in L^1 .
\]
Here $\Div$ takes the divergence of each row. For $0<h\le1$,
\[
  \mathcal G_h(\mathcal A)\le C_{J_1}
  \bigl(\|\tr\mathcal A\|_1+\|\mathcal E\|_2+\|\mathcal F\|_1\bigr)^{3/2}
  (1+|\log h|)^{3/8}.
\]
Here $C_{J_1}$ depends only on the kernel and the dimension; it is
independent of $h$ and the supports of the fields.
\end{proposition}
For determinant estimates on positive divergence-free tensors, see
Serre \cite{Serre2018}. Here $\Div\mathcal A$ includes the divergence of
an $L^2$ field, which is handled by a duality specific to the line
projection.

After compact localization, the Fourier slice identity
\eqref{app:slice} controls the projected stress in the squared average
over directions. For every compactly supported $f\in L^2(\R^4)$,
componentwise for matrix fields,
\[
 \int_{\mathbb S^3}
   \|X_\omega f\|_{\dot H^{1/2}(\omega^\perp)}^2\dd\omega
 =C\|f\|_{L^2(\R^4)}^2.
\]
In the three transverse variables, a periodic Monge--Amp\`ere test
$\psi$ has a mollified Hessian controlled in $\dot H^{-1/2}$, with
the logarithmic loss in \eqref{rt:trace-D}. The local mass growth of
$|\Delta\psi|$ is at the critical order for the kernel $|y|^{-2}$;
integrating over intermediate scales produces the logarithm. Pairing
these two Sobolev bounds controls the viscous term. Positivity then
converts the tested balance into the determinant estimate.
Appendix~\ref{app:middle-tools} proves these steps and justifies the
low-frequency estimates and removal of the localization. Only the
angular average is used; no bound uniform in each individual
direction is required.

Convolve the mass and momentum equations in space with the Poisson kernel:
\[
 \eta_r(x)=\frac{r}{\pi^2(r^2+|x|^2)^2},\qquad
 \widehat{\eta_r}(\xi)=e^{-r|\xi|},\qquad
 f_r=\eta_r*f,\quad P_r=\eta_r*P.
\]
The subscript $r$ applies to the entire expression: for example,
$(\rho u\otimes u)_r=\eta_r*(\rho u\otimes u)$ and
$(\rho|u|^2)_r=\eta_r*(\rho|u|^2)$. Use this same positive kernel
for every tensor entry. With the
cutoff $\zeta$ fixed independently of all scales and momentum levels, put
\[
 \begin{aligned}
 \mathcal A_r&=\zeta
 \begin{pmatrix}\rho_r&m_r^{\top}\\m_r&(\rho u\otimes u)_r+P_r\mathrm{Id}_3\end{pmatrix},
 \\
 \mathcal E_r&=\zeta
 \begin{pmatrix}0&0\\0&S_r\end{pmatrix}.
 \end{aligned}
\]
After extension by zero in time, positivity of $\eta_r$,
Jensen's inequality and \eqref{ao:balance} give
\[
 \mathcal A_r\ge0,\qquad P_r\ge a\rho_r^\gamma,\qquad
 \Div(\mathcal A_r-\mathcal E_r)=\zeta'\cdot(\rho_r,m_r).
\]
The three norms in Proposition~\ref{rt:optimized-CI} are bounded uniformly in
$r$ by \eqref{ao:natural} and Young's inequality. Thus
\begin{equation}\label{ao:projection}
 \mathcal G_h(\mathcal A_r)\le C_*(1+|\log h|)^{3/8},
 \qquad r>0,\quad 0<h\le1.
\end{equation}

\subsection{A cover of the momentum superlevel set}
\label{ao:cover}

Fix \(0<r<1\). Partition \(\R^3\) into cubes \(Q_j\) of side \(r\).
Choose smooth cutoffs of bounded overlap satisfying
\[
 \mathbf 1_{Q_j+B_r}\le\phi_j\le\mathbf 1_{Q_j+B_{2r}},\qquad
 f_j(t)=\frac1r\int_{r/2}^{2r}\!\int\phi_j(y)|m_s(t,y)|\dd y\dd s .
\]
For a.e. $t$, $m_s(t,x)$ is harmonic in $(x,s)$, so $|m_s(t,x)|$ is
subharmonic. The mean-value inequality and bounded overlap give, respectively,
\begin{equation}\label{ao:submean}
 |m_r(t,x)|\le Cr^{-3}f_j(t),\qquad x\in Q_j.
\end{equation}
\begin{equation}\label{ao:envelope-mass}
 \sum_j f_j(t)\le C\|m(t)\|_1,\qquad
 \sum_j\int_I f_j(t)\,\mathrm dt\le C_*.
\end{equation}
For the time level sets, always take the absolutely continuous representative
of $f_j$ on $[t_-,t_+]$ given by \eqref{ao:variation}. Thus every strict
superlevel set of $f_j$ in $I$ is relatively open, and \eqref{ao:submean} holds
outside a common null set for the countably many cubes.
Choose \(c_0\) from \eqref{ao:submean}. For each \(j\), keep the components
of \(\{f_j>c_0\Lambda r^3/2\}\cap I\) that contain a point where \(f_j\ge c_0\Lambda r^3\).
Let \(U_\Lambda\) be the union of their products with \(Q_j\), \(N\) their
number, and \(\ell\) the sum of their lengths. Then
\begin{equation}\label{ao:active-volume}
 |U_\Lambda|=r^3\ell,\qquad
 \{|m_r|\ge \Lambda\}\cap(I\times\R^3)\subset U_\Lambda
 \quad\text{up to null sets}.
\end{equation}
Each product lies in a tube of radius \(Cr\) about its vertical centre segment.

\begin{lemma}[Number and length of the intervals]\label{ao:cover-lemma}
The collection is finite and, for \(\Lambda\ge1\),
\begin{equation}\label{ao:coarse-complexity}
 \ell+rN\le \frac {C_*}{\Lambda r^3}.
\end{equation}
For \(t\in I\), set
\[
 \begin{aligned}
 w_j(x)&=\frac1r\int_{r/2}^{2r}(\eta_s*\phi_j)(x)\dd s,\\
 w(t,x)&=\sum_{j:\,\{t\}\times Q_j\subset U_\Lambda}w_j(x).
 \end{aligned}
\]
Up to cube boundaries, the weight satisfies
\begin{equation}\label{ao:raw-count}
 \ell+rN\le \ell+\frac {C_*}{\Lambda r^3}
 \left(\iint w\cdot\bigl(|\rho u\otimes u+P\mathrm{Id}_3|+|S|\bigr)+r\right),
\end{equation}
where unmarked double integrals are over \(I\times\R^3\).
Moreover,
\begin{equation}\label{ao:weight-bounds}
 0\le w\le C,\qquad \iint w\le Cr^3\ell.
\end{equation}
\end{lemma}
The proof in Appendix~\ref{app:momentum-cover} applies the momentum
equation on each selected interval.
\subsection{A lower bound on the whole cover}
\label{ao:full-union}

A cube may contain points where \(m_r=0\), even on a selected interval.
The following lower bound holds on the whole cover, up to a null set.
\begin{lemma}[Density and kinetic energy on the cover]\label{ao:active-product}
We have
\begin{equation}\label{ao:product}
 \rho_r\cdot(\rho|u|^2)_r\ge c\Lambda^2
 \qquad\text{for a.e. }(t,x)\in U_\Lambda.
\end{equation}
There is a fixed \(C_1=C_1(C_0)\ge1\) such that, for \(\Lambda\ge C_1r^{-3/2}\),
\begin{equation}\label{ao:density-lower}
 \rho_r\ge C_1^{-2}\Lambda^2r^3\ge1
 \qquad\text{for a.e. }(t,x)\in U_\Lambda.
\end{equation}
\end{lemma}
\begin{proof}
For a.e. $t\in I$, the following inequalities hold for every $x\in Q_j$,
every $y\in\supp\phi_j$, and every $s\in[r/2,2r]$:
\[
 |m_s(y)|^2\le \rho_s(y)\cdot(\rho|u|^2)_s(y),\qquad
 \eta_s(y-z)\le C\eta_r(x-z).
\]
Consequently,
\[
 \begin{aligned}
 f_j(t)&\le Cr^3\sqrt{\rho_r(t,x)\cdot(\rho|u|^2)_r(t,x)},\\
 (\rho|u|^2)_r(t,x)&\le \|\eta_r\|_\infty\|\rho(t)|u(t)|^2\|_1\le C C_0r^{-3}.
 \end{aligned}
\]
On \(U_\Lambda\), \(f_j>c_0\Lambda r^3/2\); the two claims follow.
\end{proof}

\begin{lemma}[Lower bound for the projected determinant]\label{ao:minor}
Let \(E\subset I\times\R^3\) be measurable. Suppose, for a fixed $c_*>0$,
\[
 \rho_r\ge A\ge1,\qquad \rho_r\cdot(\rho|u|^2)_r\ge c_*\Lambda^2\quad\text{a.e. on }E .
\]
Then for \(\omega=(\omega_0,\omega_x)\in\mathbb S^3\) with
\(|\omega_0|\ge1/2\),
\begin{equation}\label{ao:minor-bound}
 \det\!\left(J_h*\Pi_\omega X_\omega
                   \mathcal A_r\Pi_\omega^{\top}\right)^{1/2}
 \ge c_{a,c_*} A^{\gamma-\frac12}\Lambda
       (J_h*X_\omega\mathbf1_E)^{3/2}.
\end{equation}
The constant $c_{a,c_*}>0$ is independent of $r,h,A,\Lambda,E$ and $\omega$.
\end{lemma}
\begin{proof}
Fix $\omega$ as above and $y\in\omega^\perp$. Write $\nu$ for the
measure on $\R^4$ with density $J_h(y-\Pi_\omega(t,x))$ representing
$J_h*X_\omega$ at $y$. Put $P_E=\int_E P_r\dd\nu$, without
normalizing by $\nu(E)$. For $v=\Pi_\omega(1,0)$ the projected pressure
matrix is
\[
 Q_\omega=\Pi_\omega\operatorname{diag}(0,\mathrm{Id}_3)\Pi_\omega^{\top}
   =\mathrm{Id}_3-v\otimes v,\qquad
 \det Q_\omega=\omega_0^2,\qquad
 \operatorname{adj}Q_\omega=\omega_0^2\mathrm{Id}_3+v\otimes v .
\]
Since $|v|^2=1-\omega_0^2$, for $z=\Pi_\omega(1,u)$ we have
\begin{equation}\label{app:minor-identity}
 z^{\top}\operatorname{adj}(Q_\omega)z
   =\omega_0^2|z|^2+(v\cdot z)^2=|\omega_x-\omega_0u|^2 .
\end{equation}
On $E$ the cutoff equals one. Restrict the projected tensor integral
to $E$; it decreases in the positive semidefinite order. The resulting
matrix is $P_E Q_\omega+\Sigma$, where
\[
 \Sigma=\int_E (\rho\,z\otimes z)_r\dd\nu\ge0 .
\]
Since $Q_\omega>0$, the eigenvalues of $Q_\omega^{-1/2}\Sigma Q_\omega^{-1/2}$ are nonnegative.
Retaining the cubic and quadratic powers of $P_E$ in the determinant gives
\[
 \begin{aligned}
 \det(J_h*(\mathcal A_r)_\omega)&\ge \det(P_E Q_\omega+\Sigma)
  \\
 &\ge P_E^{3}\det Q_\omega
       +P_E^{2}\operatorname{adj}(Q_\omega):\Sigma\\
 & = P_E^{2}\int_E
 \bigl[\omega_0^2 P_r+(\rho|\omega_x-\omega_0 u|^2)_r\bigr]\dd\nu .
 \end{aligned}
\]
The velocity variance remains in the second term, even if \(m_r=0\):
\[
 (\rho|\omega_x-\omega_0u|^2)_r
 =\rho_r\left|\omega_x-\omega_0\frac{m_r}{\rho_r}\right|^2
  +\omega_0^2\left((\rho|u|^2)_r-\frac{|m_r|^2}{\rho_r}\right).
\]
Since \(|\omega_0|\ge1/2\) and \(P_r\ge a\rho_r\) on \(E\),
\[
 \begin{aligned}
 (\rho|u|^2)_r&\le\frac{2}{\omega_0^2}
       \bigl[(\rho|\omega_x-\omega_0u|^2)_r+|\omega_x|^2\rho_r\bigr]\\
 &\le C_a\bigl[(\rho|\omega_x-\omega_0u|^2)_r+\omega_0^2P_r\bigr].
 \end{aligned}
\]
Since $E\subset I\times\R^3$, $\nu(E)\le\frac{|I|}{|\omega_0|}$.
If $\nu(E)=0$, there is nothing to prove. Otherwise
$(\rho|u|^2)_r\ge c_*\Lambda^2/\rho_r$ and Jensen's inequality for
$\frac{\mathbf1_E\,\dd\nu}{\nu(E)}$ give
\[
 \begin{aligned}
 P_E^{2}\int_E (\rho|u|^2)_r\dd\nu
 &\ge c_{a,c_*}\Lambda^2\left(\int_E\rho_r^\gamma\dd\nu\right)^2
                           \int_E\rho_r^{-1}\dd\nu\\
 &\ge c_{a,c_*}\Lambda^2\nu(E)^3
       \left(\frac{1}{\nu(E)}\int_E\rho_r\dd\nu\right)^{2\gamma-1}
 \ge c_{a,c_*}\Lambda^2\nu(E)^3 A^{2\gamma-1}.
 \end{aligned}
\]
Take square roots and use $\nu(E)=J_h*X_\omega\mathbf1_E$.
\end{proof}

For $h=r$, the projection of a selected interval times a cube,
enlarged by $B_r$ in $\omega^\perp$, lies in at most
$C(1+|J'|/r)$ balls of radius $Cr$, where $J'$ is that time interval.
Summing over the cover and using Fubini's theorem gives
\[
 |\supp(J_r*X_\omega\mathbf1_{U_\Lambda})|\le Cr^2(\ell+rN),\qquad
 \int_{\omega^\perp} J_r*X_\omega\mathbf1_{U_\Lambda}=|U_\Lambda|.
\]
Thus, for $|U_\Lambda|>0$, H\"older's inequality yields
\begin{equation}\label{r2:shadow}
 \int_{\omega^\perp}(J_r*X_\omega\mathbf1_{U_\Lambda})^{3/2}
    \ge c\frac{|U_\Lambda|^{3/2}}{r(\ell+rN)^{1/2}}.
\end{equation}
Take \(A=C_1^{-2}\Lambda^2r^3\) in \eqref{ao:minor-bound}, integrate over
\(|\omega_0|\ge1/2\), and use \eqref{ao:projection}. This gives
\begin{equation}\label{ao:occupation}
 |U_\Lambda|\le C_* r^{5/3-2\gamma}\Lambda^{-4\gamma/3}
                    (\ell+rN)^{1/3}(1+|\log r|)^{1/4},
 \qquad \Lambda\ge C_1r^{-3/2}.
\end{equation}
\subsection{The flux on the selected time intervals}\label{ao:feedback}

Since \(\bigl|\{t:w(t,\cdot)\ne0\}\bigr|\le \ell\), \eqref{ao:natural}
and \eqref{ao:weight-bounds} give
\[
 \iint w|\rho u\otimes u+P\mathrm{Id}_3|\le C_* \ell,
 \qquad
 \iint w|S|\le C_* r^{3/2}\ell^{1/2}.
\]
Substitution in \eqref{ao:raw-count} yields
\begin{equation}\label{ao:active-count}
 \ell+rN\le \ell+\frac {C_*}{\Lambda r^3}
                   \left(\ell+r^{3/2}\ell^{1/2}+r\right).
\end{equation}
Young's inequality and \eqref{ao:envelope-mass} give, respectively,
\[
 \|m_r(t)\|_\infty\le C_* r^{-\frac{3(\gamma+1)}{2\gamma}},\qquad
 U_\Lambda\ne\varnothing\ \Longrightarrow\ \Lambda r^3\le C_*.
\]
Inserting \(\ell=r^{-3}|U_\Lambda|\) in \eqref{ao:active-count} yields
\begin{equation}\label{ao:complexity-feedback}
 \ell+rN\le C_*\left(
       \frac{|U_\Lambda|}{\Lambda r^6}
       +\frac{|U_\Lambda|^{1/2}}{\Lambda r^3}
       +\frac1{\Lambda r^2}\right).
\end{equation}
Cube \eqref{ao:occupation} and use this bound:
\[
 |U_\Lambda|^3\le C_*(1+|\log r|)^{3/4}\Lambda^{-4\gamma-1}
 \left(r^{-1-6\gamma}|U_\Lambda|+r^{2-6\gamma}|U_\Lambda|^{1/2}+r^{3-6\gamma}\right).
\]
Bounding the sum by three times its largest term gives
\begin{align}
 |U_\Lambda|\le C_*\bigl[&
 r^{-1/2-3\gamma}\Lambda^{-(4\gamma+1)/2}
                         (1+|\log r|)^{3/8}\notag\\
 &+r^{(4-12\gamma)/5}\Lambda^{-(8\gamma+2)/5}
                         (1+|\log r|)^{3/10}\notag\\
 &+r^{1-2\gamma}\Lambda^{-(4\gamma+1)/3}
                         (1+|\log r|)^{1/4}\bigr].
 \label{ao:three-bounds}
\end{align}
The three terms correspond to the flux, the viscous stress and the two
endpoints of \(I\). Alternatively, \eqref{ao:occupation} and
\eqref{ao:coarse-complexity} give
\begin{equation}\label{ao:coarse-occupation}
 |U_\Lambda|\le C_* r^{2/3-2\gamma}\Lambda^{-(4\gamma+1)/3}
                              (1+|\log r|)^{1/4}.
\end{equation}
\subsection{Summation over momentum levels}
\label{ao:summation}

\begin{proof}[Proof of Proposition~\ref{ao:gain}]
Set \(\Lambda_0=C_1r^{-3/2}\). The low levels satisfy
\[
 r^2\iint_{\{|m_r|\le \Lambda_0\}}|m_r|^2
 \le r^2\Lambda_0\|m_r\|_{L^1(I\times\R^3)}
 \le C_*r^{1/2}.
\]
For nonempty dyadic shells with \(\Lambda=2^j\Lambda_0\), \(j\ge0\), use
\[
 r^2\iint_{\{\Lambda\le |m_r|<2\Lambda\}} |m_r|^2\le4r^2\Lambda^2|U_\Lambda|,\qquad
 \Lambda\le \Lambda_{\max}:=C_*r^{-\frac{3(\gamma+1)}{2\gamma}}.
\]
If \(9/8<\gamma<5/4\), the coarse bound \eqref{ao:coarse-occupation}
and the flux term in \eqref{ao:three-bounds}, after multiplication by
\(r^2\Lambda^2\) and omission of logarithms, are
\[
 r^{\frac83-2\gamma}\Lambda^{\frac{5-4\gamma}{3}},\qquad
 r^{\frac32-3\gamma}\Lambda^{\frac{3-4\gamma}{2}} .
\]
The first increases in \(\Lambda\), the second decreases. They agree at
\[
 \Lambda_*=r^{-\frac{7+6\gamma}{4\gamma+1}},\qquad
 \Lambda_0<\Lambda_*<\Lambda_{\max}\quad(r\text{ small}),
\]
since
\[
 \begin{aligned}
 \frac{7+6\gamma}{4\gamma+1}-\frac32
    &=\frac{11}{2(4\gamma+1)}>0,\\
 \frac{3(\gamma+1)}{2\gamma}-\frac{7+6\gamma}{4\gamma+1}
    &=\frac{\gamma+3}{2\gamma(4\gamma+1)}>0.
 \end{aligned}
\]
Sum the first below \(\Lambda_*\), the second above \(\Lambda_*\).
Both geometric series have the power
\begin{equation}\label{ao:critical-power}
 \frac83-2\gamma
 -\frac{7+6\gamma}{4\gamma+1}\frac{5-4\gamma}{3}
 =\frac{8\gamma-9}{4\gamma+1}.
\end{equation}
The two remaining terms in \eqref{ao:three-bounds} may be summed over
all levels since \(\min(a,b+c+d)\le\min(a,b)+c+d\).
Their powers of \(\Lambda\) have opposite signs:
\[
 \begin{aligned}
 r^{(14-12\gamma)/5} \Lambda_0^{(8-8\gamma)/5}
   &\le C_* r^{2/5} &&\text{(stress)},\\
 r^{3-2\gamma}\Lambda_{\max}^{(5-4\gamma)/3}
   &\le C_* r^{\frac{5(\gamma-1)}{2\gamma}} &&\text{(endpoints)} .
 \end{aligned}
\]
Both powers of \(r\) exceed \eqref{ao:critical-power}. The logarithmic
factor for each sum is at most \((1+|\log r|)^{3/8}\).

For \(\gamma>5/4\), \eqref{ao:coarse-occupation} alone gives a decreasing
geometric sum of order \(r^{1/6}(1+|\log r|)^{1/4}\).
At \(\gamma=5/4\) its terms are constant in \(\Lambda\), with at most
\(C_*(1+|\log r|)\) levels. Therefore for every \(\gamma>9/8\),
\begin{equation}\label{ao:scale-bound}
 r^2\|m_r\|_{L^2(I\times\R^3)}^2
 \le C_*\left[
 r^{1/2}
 +r^{\min\{\frac16,\frac{8\gamma-9}{4\gamma+1}\}}
                    (1+|\log r|)^{5/4}\right].
\end{equation}
Young's inequality covers \(r\) bounded away from zero. Plancherel gives
\[
 \int_0^\infty r^{2-2\alpha}
       \|\eta_r*m\|_{L^2(I\times\R^3)}^2\,\frac{\mathrm dr}{r}
 =\widetilde C_\alpha\|I_{1-\alpha}m\|_{L^2(I\times\R^3)}^2,
 \qquad 0<\widetilde C_\alpha<\infty ,
\]
since
\[
 \int_0^\infty r^{1-2\alpha}e^{-2r|\xi|}\,\mathrm dr
 =|\xi|^{-2+2\alpha}\int_0^\infty s^{1-2\alpha}e^{-2s}\,\mathrm ds .
\]
At \(r=0\), \eqref{ao:scale-bound} is integrable against
\(r^{-2\alpha}\dd r/r\) for the indicated range of \(\alpha\).
For \(r\ge1\), Young's inequality gives \(\|m_r(t)\|_2^2\le C_*r^{-3}\).
Lemma~\ref{kt:tails} below gives the two kinetic energy estimates.
\end{proof}

\subsection{Kinetic energy at large velocities and densities}\label{kt:section}

The truncation $(1-M/|u|)_+\cdot u$ is parallel to $m$. This gives a
nonnegative test for the kinetic energy on $\{|u|>2M\}$.

\begin{lemma}[Kinetic energy on superlevel sets]\label{kt:tails}
Let \(I\) be bounded, \(\gamma>1\), and \(0<\alpha<1\). Suppose
\[
 \rho\ge0,\qquad
 \rho\in L^\infty(I;L^1(\R^3)\cap L^\gamma(\R^3)),\qquad
 \rho|u|^2\in L^\infty(I;L^1(\R^3)),
\]
\[
 u\in L^2(I;\dot H^1(\R^3)),\qquad
 I_{1-\alpha}(\rho u)\in L^2(I\times\R^3).
\]
Use the \(L^6\) representative of \(u(t)\). Then for \(M,A\ge1\),
\begin{align}
 \iint_{\substack{t\in I\\|u|>M}}\rho|u|^2
 &\le C_{*,\alpha} M^{-\frac{2\alpha}{1+\alpha}},\label{kt:velocity-tail}\\
 \iint_{\substack{t\in I\\\rho>A}}\rho|u|^2
 &\le C_{*,\alpha} A^{-\frac{(\gamma-1)\alpha}{1+2\alpha}}.\label{kt:density-tail}
\end{align}
Here and in the following corollary, $C_{*,\alpha}$ depends only on
$I,\gamma,\alpha$ and the indicated norms, and is independent of the
levels $M,A,L$ and the truncation index $N$. No equation of motion is needed.
\end{lemma}

\begin{proof}
Let
\[
 \Phi_M(v)=\left(1-\frac M{|v|}\right)_+\cdot v,\qquad \Phi_M(0)=0.
\]
As \(\Phi_M\) is 1-Lipschitz, the chain rule and Sobolev inequality give
\[
\begin{aligned}
 \|\Phi_M(u)\|_2^2
 &\le \int_{\{|u|>M\}} |u|^2
 \le M^{-4}\|u\|_6^6\le C M^{-4}\|\nabla u\|_2^6,\\
 \|\nabla\Phi_M(u)\|_2&\le \|\nabla u\|_2 .
\end{aligned}
\]
Interpolation between \(L^2\) and \(\dot H^1\) gives
\begin{equation}\label{kt:interpolation}
 \|\Phi_M(u(t))\|_{\dot H^{1-\alpha}}
 \le C M^{-2\alpha}\|\nabla u(t)\|_2^{1+2\alpha}.
\end{equation}
For a.e. \(t\), let \(v_N\) be the radial truncation of \(\Phi_M(u(t))\)
at length \(N\). Then
\[
 v_N\to\Phi_M(u)\quad\text{in }H^1\text{ and in }\dot H^{1-\alpha},
 \qquad m\cdot v_N\uparrow\rho|u|^2\cdot(1-M/|u|)_+ .
\]
Cutoff and mollification approximate \(v_N\) in \(H^1\) and a.e. by
smooth fields of size at most \(N\). Since \(m(t)\in L^1\), dominated
convergence identifies \(\langle m,v_N\rangle=\int m\cdot v_N\).
Pass to \(N\to\infty\) by duality and monotone convergence:
\[
 \frac12\int_{\{|u|>2M\}}\rho|u|^2
 \le \langle m,\Phi_M(u)\rangle
 \le \|I_{1-\alpha}m\|_2
      \|\Phi_M(u)\|_{\dot H^{1-\alpha}}.
\]
For \(L\ge1\), Cauchy--Schwarz and \eqref{kt:interpolation} give
\[
 \int_{\{\|\nabla u(t)\|_2^2\le L\}}\int_{\{|u|>2M\}}\rho|u|^2
 \le C M^{-2\alpha}L^\alpha
       \|I_{1-\alpha}m\|_{L^2(I\times\R^3)}
       \|\nabla u\|_{L^2(I\times\R^3)}.
\]
On the remaining times, Chebyshev's inequality and the energy bound give
\[
 \int_{\{\|\nabla u(t)\|_2^2>L\}}\!\!\int\rho|u|^2
 \le \frac{1}{L}\|\nabla u\|_{L^2(I\times\R^3)}^2
                  \|\rho|u|^2\|_{L^\infty(I;L^1)} .
\]
Hence
\[
 \iint_{\substack{t\in I\\|u|>2M}}\rho|u|^2
 \le C_{*,\alpha}\bigl(M^{-2\alpha}L^\alpha+L^{-1}\bigr).
\]
Take \(L=M^{2\alpha/(1+\alpha)}\) and use the energy bound for thresholds
in \([1,2]\). This proves \eqref{kt:velocity-tail}. For \eqref{kt:density-tail},
\[
 \iint_{\substack{t\in I\\\rho>A}}\rho|u|^2
 \le M^2\iint_{\substack{t\in I\\\rho>A}}\rho
      +\iint_{\substack{t\in I\\|u|>M}}\rho|u|^2
 \le C_{*,\alpha}\left(M^2A^{1-\gamma}+M^{-\frac{2\alpha}{1+\alpha}}\right).
\]
The choice
\[
 M=A^{\frac{(\gamma-1)(1+\alpha)}{2(1+2\alpha)}}
\]
makes the two powers equal and proves \eqref{kt:density-tail}.
\end{proof}

\begin{corollary}[Higher integrability of the kinetic energy density]
\label{kt:kinetic-integrability}
Under the hypotheses of Lemma~\ref{kt:tails},
\[
 \iint_{\{\rho|u|^2>L\}}\rho|u|^2
 \le C_{*,\alpha} L^{-\frac{(\gamma-1)\alpha}{\gamma+(\gamma+1)\alpha}},
 \qquad L\ge1.
\]
In particular, \(\rho|u|^2\in L^{1+s}(I\times\R^3)\) for
\[
 0<s<\frac{(\gamma-1)\alpha}{\gamma+(\gamma+1)\alpha}.
\]
The constant in the superlevel estimate has the dependence specified in
Lemma~\ref{kt:tails}. The $L^{1+s}$ bound has a constant $C_{*,\alpha,s}$,
which may additionally depend on $s$. For each fixed admissible $s$,
the bounds are uniform for families satisfying the hypotheses uniformly.
\end{corollary}
\begin{proof}
For \(1\le A\le L\),
\[
 \{\rho|u|^2>L\}\subset
 \{\rho>A\}\cup\{|u|>(L/A)^{1/2}\}.
\]
Apply \eqref{kt:velocity-tail} and \eqref{kt:density-tail} and take
\[
 A=L^{\frac{1+2\alpha}{\gamma+(\gamma+1)\alpha}} .
\]
This proves the superlevel estimate. The layer-cake formula gives
\[
 \iint (\rho|u|^2)^{1+s}
 =s\int_0^\infty L^{s-1}
             \iint_{\{\rho|u|^2>L\}}\rho|u|^2\,\mathrm dL .
\]
Use the \(L^1\) bound on \(0<L<1\) and the superlevel estimate on \(L\ge1\).
\end{proof}

\subsection{The endpoint in this argument}\label{ao:threshold}

At \(\gamma=9/8\), \eqref{ao:critical-power} vanishes and
\(\Lambda_*=r^{-5/2}\). For a fixed smooth spatial profile rescaled to a ball
of radius \(r\), the corresponding orders are
\[
 \rho\sim r^{-2},\qquad |u|\sim r^{-1/2},\qquad
 \rho|u|^2\sim r^{-3},\qquad |m|\sim r^{-5/2}.
\]
Thus in this scaling
\[
 \int \rho|u|^2\sim1,\qquad \int|\nabla u|^2\sim1,\qquad
 \|m\|_{\dot H^{-1}}\sim1,\qquad \int\rho^\gamma\sim r^{3-2\gamma}.
\]
These rescaled profiles are not asserted to solve the equations.
The vanishing exponent is a limitation of the two estimates in
\eqref{ao:critical-power}; it gives no nonexistence result at the endpoint.

For \(\gamma>9/8\), Proposition~\ref{ao:gain} supplies
\eqref{cp:momentum_gain} and the density tail \eqref{r6:assumed-kinetic-tail}
for the whole-space approximations. Theorem~\ref{cp:compactness} identifies
the limit; Sections~\ref{r4:global-energy} and~\ref{sec:construction}
pass the energy inequality to that limit and complete
Theorem~\ref{thm:nine-eighths}.

\section{Pressure estimates and energy convergence on \texorpdfstring{$\R^3$}{R3}}
\label{r4:global-energy}

We pass to the energy inequality between positive times by proving
\[
 E_j\longrightarrow E_\delta\quad\text{in }L^1(I),
 \qquad I\Subset(0,\infty).
\]
A pressure moment and an estimate on the kinetic energy at high density
control the energy outside $B_R$. All spatial norms are over $\R^3$.

\subsection{Higher integrability of the pressure}

\begin{lemma}[Commutator estimate on the whole space]
\label{r4:homogeneous-commutator}
Let $3<p<\infty$, $s=1-3/p$, and $R_{ik}=\partial_i\partial_k\Delta^{-1}$.
If $v\in \dot H^1(\R^3)\cap L^6(\R^3)$ and $f\in L^p\cap L^3$, then
\begin{equation}\label{r4:homogeneous-commutator-bound}
 \|[v,R_{ik}]f\|_{H^s}
 \le C_p\|\nabla v\|_2\bigl(\|f\|_p+\|f\|_3\bigr).
\end{equation}
\end{lemma}

\begin{proof}
For $v\in C_c^\infty$, apply Lemma~\ref{cp:commutator} to
$v(L\,\cdot)$ and $f(L\,\cdot)$. Since $R_{ik}$ commutes with dilations
and $s=1-3/p$, division by $L^{s-3/2}$ gives
\[
 \|[v,R_{ik}]f\|_{\dot H^s}
 \le C_p\bigl(L^{-1}\|v\|_2+\|\nabla v\|_2\bigr)\|f\|_p.
\]
Let $L\to\infty$. For the $L^2$ norm, the Riesz transform bounds and
Sobolev's inequality give
\[
 \|vR_{ik}f-R_{ik}(vf)\|_2
 \le C\|v\|_6\|f\|_3
 \le C\|\nabla v\|_2\|f\|_3.
\]
Approximation in $\dot H^1$ gives the result for general $v$. The
second estimate identifies the $L^2$ limit with the stated commutator.
\end{proof}

\begin{proposition}[Pressure estimate on interior time intervals]
\label{r4:global-pressure}
Let $(\rho_j,u_j)$ satisfy the momentum equation and the renormalized
continuity equation on $\R^3$, with
\[
 P_j=a\rho_j^\gamma+\delta_j\rho_j^\beta,\qquad
 \gamma>1,\quad \beta>\max\{4,\gamma\},\quad \delta_j>0.
\]
Assume uniform bounds for the mass, pressure, kinetic energy and
viscous dissipation, and assume $\rho_j\in L^\infty_tL^\beta_x$
for each $j$. For each $j$, $u_j(t)$ is taken to be the
$L^6(\R^3)$ representative of its class in $\dot H^1(\R^3)$
for almost every $t$. Suppose that for some $0<\sigma<1$ and every interior
time interval $J\Subset(0,\infty)$,
\[
 \|m_j\|_{L^2(J;H^{-1+\sigma}(\R^3))}\le C_{J,\sigma} .
\]
For every $I\Subset J$,
\begin{equation}\label{r4:global-pressure-moment}
 \int_I\int_{\R^3}P_j\rho_j^\theta\,\dd x\,\dd t
 \le C_{I,J,\sigma,\theta},
 \qquad 0<\theta<1,\quad \theta\le\frac{\gamma\sigma}{3}.
\end{equation}
The constant may also depend on the fixed physical parameters, $\gamma,\beta$,
and the stated uniform bounds, including the momentum bound on $J$.
It is independent of $j$, and is uniform in the coefficients $\delta_j$
whenever those stated bounds are uniform.
\end{proposition}

\begin{proof}
Set $p=3/\sigma$. Then $p>3$, $\theta p\le\gamma$ and $1-3/p=1-\sigma$.
In the estimates below, $C_*$ may depend on the data listed in the proposition
and the fixed cutoffs, but not on $j,N,R$.
Use the bounded increasing approximations $b_N$ of
$(1+z)^\theta-1$ from the proof of Theorem~\ref{cp:compactness}. Set
\[
 f_j=b_N(\rho_j),\qquad
 g_j=\bigl(b_N(\rho_j)-\rho_j b_N'(\rho_j)\bigr)\diver u_j .
\]
The renormalized equation gives
\[
 \partial_t f_j+\diver(f_j u_j)=g_j .
\]
The mass and pressure bounds, together with
\[
 |b_N(z)|+|b_N(z)-zb_N'(z)|\le 2\min\{z,z^\theta\},
\]
give, uniformly in $j,N$,
\[
 \begin{aligned}
 \|f_j\|_{L^\infty(J;L^1\cap L^p)}
 +\|b_N(\rho_j)-\rho_j b_N'(\rho_j)\|_{L^\infty(J;L^1\cap L^p)}
 &\le C_*,\\
 \|g_j\|_{L^2(J;L^1\cap L^{2p/(p+2)})}&\le C_*.
 \end{aligned}
\]
By Lemmas~\ref{cp:potential} and \ref{r4:homogeneous-commutator},
\[
 \|\nabla\Delta^{-1}g_j\|_{L^2(J;H^{1-3/p})}
 +\|[u_{j,k},R_{ik}]f_j\|_{L^2(J;H^{1-3/p})}\le C_* .
\]
The low frequencies in Lemma~\ref{cp:potential} use only $\|g_j\|_1$.
Splitting $|x|^{-2}$ at $|x|=1$ and $|\xi|^{-1}$ at $|\xi|=1$ gives
\[
 \|\nabla\Delta^{-1}f_j\|_{L^\infty(J;L^\infty\cap L^2)}
 \le C_p\|f_j\|_{L^\infty(J;L^1\cap L^p)}.
\]

Take $0\le\zeta\in C_c^\infty(J)$ with $\zeta=1$ on $I$.
Let $\chi_R(x)=\chi(x/R)$, where $\chi$ is radial, equals one on
$B_1$, vanishes outside $B_2$, and decreases with the radius. Test
the momentum equation with
\[
 \zeta\chi_R\nabla\Delta^{-1}f_j .
\]
Here $f_j$ need not have compact support. For fixed $j,N$, the mass,
the artificial pressure, and $f_j\in L^\infty\cap L^\infty_tL^1_x$ give
\[
 \begin{gathered}
 m_j,\ f_j u_j\in L^2L^2,\qquad g_j\in L^2(L^1\cap L^2),\qquad
 \partial_t\nabla\Delta^{-1}f_j\in L^2L^2,\\
 \rho_j u_j\otimes u_j\in L^1L^r,\qquad
 \frac1r=\frac1\beta+\frac13<1 .
 \end{gathered}
\]
To justify the test, mollify $f_j$ in space and time on a slightly
larger interior interval, keeping $j,N,R$ fixed. The time derivatives of
the potentials converge in $L^2L^2$, and the spatial derivatives
$R_{ik}f_j$ of these potentials converge in $L^{r'}$ at a.e. time along a subsequence,
boundedly in $L^\infty_tL^{r'}_x$. This gives the convective pairing.
The mollified $f_j$ are bounded by $\|b_N\|_\infty$ and converge a.e.,
so dominated convergence applies to the principal pressure term.
For the term with $\nabla\chi_R$, the potentials converge locally in
$L^2_tL^\infty_x$ by the kernel estimate above, and $P_j\in L^\infty_tL^1_x$.
The inertial terms combine into the pairing
\[
 \sum_i\left\langle\zeta\chi_R m_{j,i},\,
     \partial_i\Delta^{-1}g_j+
      \sum_k[u_{j,k},R_{ik}]f_j\right\rangle .
\]
The second factor is bounded in $L^2(J;H^{1-\sigma})$ since $1-3/p=1-\sigma$.
It therefore pairs with $m_j$ in $L^2(J;H^{-1+\sigma})$; multiplication
by $\chi_R$ is bounded on $H^{1-\sigma}$ uniformly for $R\ge1$.
The stress and time cutoff contribute at most
\[
 C_*\|S(u_j)\|_{L^2L^2}\|f_j\|_{L^2L^2}
 +C_*\|m_j\|_{L^1L^1}\|\nabla\Delta^{-1}f_j\|_{L^\infty L^\infty}.
\]
The terms with $\nabla\chi_R$ are bounded by
\[
 \frac{C_*}{R}\left(\|\rho_j|u_j|^2\|_{L^1L^1}
 +\|P_j\|_{L^1L^1}
 +\|S(u_j)\|_{L^2L^2}\|\nabla\Delta^{-1}f_j\|_{L^2L^2}\right).
\]
Thus $\iint\zeta\chi_RP_j b_N(\rho_j)\le C_{I,J,\sigma,\theta}$. Sending first
$R\to\infty$, then $N\to\infty$, gives
\[
 \iint_{I\times\R^3} P_j\cdot\big[(1+\rho_j)^\theta-1\big]\le C_{I,J,\sigma,\theta} .
\]
Since $(1+z)^\theta-1\ge c_\theta z^\theta$ for $z\ge1$ and
$P_j$ is uniformly bounded in $L^1(I\times\R^3)$, this proves
\eqref{r4:global-pressure-moment}.
\end{proof}

\subsection{Tightness at spatial infinity}

\begin{proposition}[Global convergence of the energy]
\label{r4:energy-tightness}
Assume the uniform mass and energy bounds of \eqref{eq:energy-spaces}, with the $L^6$ representatives of the velocities, and the continuity equation
on $\R^3$. Densities and velocities on expanding balls may be extended
by zero for this purpose. Let $\beta>\gamma>1$ and
$\delta_j\ge0$ with $\delta_j\to\delta\ge0$. Write
\[
 P_j=a\rho_j^\gamma+\delta_j\rho_j^\beta,\qquad
 \Theta_j(z)=\frac{a}{\gamma-1}z^\gamma+\frac{\delta_j}{\beta-1}z^\beta .
\]
Assume strong convergence of the initial densities in $L^1$, weak
convergence of the velocity gradients in global $L^2$ on finite time
intervals, and the following local convergences:
\[
 \begin{aligned}
 \rho_j&\to\rho &&\text{a.e.},&
 m_j&\to\rho u &&\text{in distributions},\\
 P_j&\to a\rho^\gamma+\delta\rho^\beta &&\text{in }L^1_{\rm loc},&
 \rho_j|u_j|^2&\rightharpoonup\rho|u|^2 &&\text{weakly in }L^1_{\rm loc}.
 \end{aligned}
\]
Suppose that the left-hand side of \eqref{r4:global-pressure-moment}
is bounded uniformly in $j$ for some $\theta>0$ on every interior
interval $I$, and that, for some $\kappa>0$ and every such $I$,
\begin{equation}\label{r6:assumed-kinetic-tail}
 \sup_j\int_I\int_{\{\rho_j>A\}}\rho_j|u_j|^2\,\dd x\,\dd t
 \le C_{I,\kappa} A^{-\kappa},\qquad A\ge2 .
\end{equation}
Then
\begin{align}
 \rho_j&\to\rho
   &&\text{in }L^1(I\times\R^3)\cap L^\gamma(I\times\R^3),
       \label{r4:global-density-convergence}\\
 a\rho_j^\gamma&\to a\rho^\gamma,\qquad
 \delta_j\rho_j^\beta\to\delta\rho^\beta
   &&\text{in }L^1(I\times\R^3),
       \label{r4:global-pressure-convergence}\\
 \sqrt{\rho_j}\,u_j&\to\sqrt\rho\,u
   &&\text{in }L^2(I\times\R^3),
       \label{r4:global-kinetic-convergence}\\
 E_j&\to E_\delta
   &&\text{in }L^1(I).
       \label{r4:global-energy-convergence}
\end{align}
Here $E_j$ is the total energy with internal energy $\Theta_j$ and $E_\delta$
is defined with coefficient $\delta$. If $\delta>0$, then
$\rho_j\to\rho$ also in $L^\beta(I\times\R^3)$.
\end{proposition}

\begin{proof}
The decorated constants in this proof may depend on the fixed pressure
parameters, the stated uniform bounds, $\sup_j\delta_j$ and a fixed
cutoff profile, but not on $j,R,A$. We write $C_*$ for dependence on
these data alone; additional time and exponent dependence is displayed
in subscripts.
\emph{Mass.} Fix $T>\sup I$. Test continuity with a cutoff of the
exterior of $B_R$, using $\|m_j\|_{L^\infty_tL^1_x}\le C_*$. This gives
\[
 \sup_j\esssup_{0<t<T}\int_{|x|>2R}\rho_j(t,x)\,\dd x
 \le\sup_j\int_{|x|>R}\rho_{0j}\,\dd x+\frac{C_{*,T}}{R}
 \longrightarrow0 .
\]
Here $j$ ranges over solutions defined up to $T$. The initial tails
vanish uniformly since $\rho_{0j}\to\rho_0$ in $L^1$. Hence
\[
 M_R:=\sup_j\esssup_{0<t<T}\int_{|x|>R}\rho_j(t)\,\dd x
 \longrightarrow0;
\]
Fatou's lemma gives the same exterior bound for $\rho$.

\emph{Pressure.} On $\{\rho_j>A\}$ use
\eqref{r4:global-pressure-moment}; on $\{\rho_j\le A\}$ use
$\rho_j^q\le A^{q-1}\rho_j$. Thus
\begin{equation}\label{r4:pressure-space-tail}
 \int_I\int_{|x|>R}P_j
 \le C_{I,\theta}A^{-\theta}
      +|I|\left(aA^{\gamma-1}+\bigl(\sup_j\delta_j\bigr)A^{\beta-1}\right)M_R .
\end{equation}
First let $R\to\infty$, then $A\to\infty$. On bounded cylinders
both pressure components converge a.e. and are uniformly integrable
by the local $L^1$ convergence of $P_j$. Vitali's theorem and
\eqref{r4:pressure-space-tail} give \eqref{r4:global-pressure-convergence}.
In particular, for $\delta=0$,
\[
 \iint_{I\times\R^3}\delta_j\rho_j^\beta
 \le C_{I,\theta}A^{-\theta}+\delta_j A^{\beta-1}\iint_{I\times\R^3}\rho_j
 \longrightarrow0\qquad(j\to\infty,\ \text{then }A\to\infty).
\]
Now $|z-w|^q\le |z^q-w^q|$ for $z,w\ge0$, $q\ge1$. This gives
$L^\gamma$ convergence and, if $\delta>0$, $L^\beta$ convergence.
Local $L^1$ convergence and $M_R\to0$ give the global $L^1$ limit.

\emph{Kinetic energy.} On the part of the exterior where $\rho_j\le A$,
\[
 \int_{\{|x|>R,\ \rho_j\le A\}}\rho_j|u_j|^2
 \le \|\rho_j\mathbf1_{\{|x|>R,\ \rho_j\le A\}}\|_{3/2}\|u_j\|_6^2
 \le A^{1/3} M_R^{2/3}\|u_j\|_6^2 .
\]
Integrate in time and use \eqref{r6:assumed-kinetic-tail} on the complement:
\begin{equation}\label{r4:kinetic-space-tail}
 \int_I\int_{|x|>R}\rho_j|u_j|^2
 \le C_{I,\kappa} A^{-\kappa}+C_{*,T} A^{1/3}M_R^{2/3}.
\end{equation}
Local weak convergence and \eqref{r4:kinetic-space-tail} give
\[
 \rho_j|u_j|^2\rightharpoonup\rho|u|^2
       \quad\text{in }L^1(I\times\R^3).
\]
Indeed, one tests first on $I\times B_R$ and then lets $R\to\infty$.
The same tail bound for the limit follows by exhaustion.

Let $\sqrt{\rho_j}\,u_j\rightharpoonup q$ in $L^2(I\times\R^3)$
along a subsequence. Then
\[
 \begin{gathered}
 \|\sqrt{\rho_j}-\sqrt\rho\|_2^2\le\|\rho_j-\rho\|_1\to0,\qquad
 m_j\to\rho u\ \text{in distributions},\\
 \sqrt\rho\,q=\rho u,\qquad
 \iint\rho|u|^2\le\|q\|_2^2\le
 \liminf_j\iint\rho_j|u_j|^2=\iint\rho|u|^2 .
 \end{gathered}
\]
Thus $q=0$ on $\{\rho=0\}$ and the norms converge. This proves
\eqref{r4:global-kinetic-convergence}. Squaring in $L^2$ and using
\eqref{r4:global-pressure-convergence} yields \eqref{r4:global-energy-convergence}.
\end{proof}

\begin{proposition}[Energy inequality between positive times]
\label{r4:differential-energy}
Under the hypotheses of Proposition~\ref{r4:energy-tightness},
suppose that each approximate solution satisfies the differential
energy inequality. Then
\begin{equation}\label{r4:differential-energy-inequality}
 -\int_0^\infty E_\delta(t)\varphi'(t)\,\dd t
 +\int_0^\infty \mathcal D(u(t))\varphi(t)\,\dd t\le0,
 \qquad 0\le\varphi\in C_c^\infty(0,\infty).
\end{equation}
In particular, for almost every pair $0<s<t$,
\[
 E_\delta(t)+\int_s^t\mathcal D(u(\tau))\,\dd\tau\le E_\delta(s).
\]
If the limit has the prescribed weak initial traces and satisfies the
energy inequality from time zero, then $E_\delta(t)\to E_{\delta,0}$
as $t\downarrow0$ through the full-measure set where that inequality holds.
\end{proposition}

\begin{proof}
For $0\le\varphi\in C_c^\infty(0,\infty)$, strong energy convergence
and weak lower semicontinuity of the quadratic dissipation give
\[
 \int E_j\varphi'\to\int E_\delta\varphi',\qquad
 \int\mathcal D(u)\varphi\le\liminf_j\int\mathcal D(u_j)\varphi.
\]
Hence \eqref{r4:differential-energy-inequality} holds and
$t\mapsto E_\delta(t)+\int_0^t\mathcal D(u)\,\dd\tau$ has a nonincreasing
representative. At Lebesgue points this gives the asserted inequality.
At the initial time, convexity, the weak traces and the energy
inequality from zero give, along the stated full-measure set,
\[
 E_{\delta,0}\le\liminf_{t\downarrow0}E_\delta(t)
 \le\limsup_{t\downarrow0}E_\delta(t)\le E_{\delta,0}.
\]
\end{proof}

\section{Construction of weak solutions}\label{sec:construction}

We first let $R\to\infty$ at fixed $\delta>0$ and then let
$\delta\downarrow0$:
\[
 p_\delta(z)=az^\gamma+\delta z^\beta,\qquad
 \beta>\max\{4,\gamma\}.
\]
At fixed $\delta$ the artificial pressure gives $\rho u\in L^2_tL^2_x$.
For the second limit we use the uniform bound of
Proposition~\ref{ao:gain}.
Constants with a $\delta$ subscript may also depend on the fixed
physical parameters, $\gamma,\beta$, and the prescribed initial mass and energy, or
the uniform bounds stated in the relevant lemma. They are uniform in
the exhaustion radius $R$, but need not be uniform as $\delta\downarrow0$.
The uniform estimates used in the second limit instead depend on the
common initial bounds and are independent of $j$.
\subsection{Expanding balls at fixed artificial pressure}

\begin{lemma}[Weak solutions at fixed artificial pressure]
\label{r4:fixed-pressure-energy}
Let $\delta>0$, $\gamma>1$ and $\beta>\max\{4,\gamma\}$. For every
nonnegative initial density and compatible momentum of finite mass
and energy for
\[
 p_\delta(z)=az^\gamma+\delta z^\beta,
\]
there is a weak solution on $\R^3$ with this pressure. Its energy
\[
 E_\delta(t)=\int_{\R^3}
  \left(\frac12\rho|u|^2+\frac{a}{\gamma-1}\rho^\gamma
                    +\frac{\delta}{\beta-1}\rho^\beta\right)(t)\,\dd x
\]
satisfies both the energy inequality from time zero and the
differential energy inequality. In this lemma the estimates may
depend on $\delta$, the fixed physical parameters, $\gamma,\beta$, the initial mass
and energy, and the finite time interval under consideration.
\end{lemma}

\begin{lemma}[Uniform pressure estimate on balls]\label{lem:ball-pressure}
Fix $\delta>0$, $\gamma>1$, $\beta>\max\{4,\gamma\}$ and $T<\infty$.
Let $(\rho_R,u_R)$ be renormalized weak solutions on $(0,T)\times B_R$
with pressure $p_\delta$ and $u_R|_{\partial B_R}=0$. Suppose
\[
 \sup_{R\ge1}\left[
 \esssup_{0<t<T}\int_{B_R}(\rho_R+\rho_R^\beta+\rho_R|u_R|^2)\dd x
 +\int_0^T\int_{B_R}|\nabla u_R|^2\dd x\dd t\right]<\infty .
\]
Then for each $I\Subset(0,T)$,
\begin{equation}\label{r4:fixed-pressure-moment}
 \int_I\int_{B_R}p_\delta(\rho_R)\rho_R\dd x\dd t\le C_{\delta,I,T},
 \qquad R\ge1 .
\end{equation}
The constant also depends on the fixed physical parameters,
$\gamma,\beta$, and the bound in the hypothesis. Its time dependence
includes the separation of $I$ from the endpoints $0,T$; it is
independent of $R$.
\end{lemma}
The Bogovski\u{\i} test proving this estimate is given in
Appendix~\ref{app:bogovskii}.

\begin{proof}[Proof of Lemma~\ref{r4:fixed-pressure-energy}]
\emph{Solutions on balls.}
On $B_R$ take the initial data restricted to the ball and impose
$u_R=0$ on $\partial B_R$. The bounded-domain existence theorem
\cite[Theorem~3.1]{Feireisl2002} applies with growth exponent $\beta$, since for suitable fixed constants
$C_{\delta,1},C_{\delta,2}>0$,
\begin{equation}\label{eq:classical-pressure-check}
 C_{\delta,1}^{-1}z^{\beta-1}-C_{\delta,2}\le p_\delta'(z)
 \le C_{\delta,1}z^{\beta-1}+C_{\delta,2},\qquad z\ge0 .
\end{equation}
Indeed $p_\delta'(z)\ge\delta\beta z^{\beta-1}$,
$z^{\gamma-1}\le1+z^{\beta-1}$ and $\mu+\lambda\ge\mu/3$.
To obtain both energy inequalities for the same solution, use
\cite[Proposition~4.1 and Definition~2.1]{FanelliFeireisl2020} with
zero force and boundary velocity, so there is no inflow. Its bulk
viscosity is $\lambda+2\mu/3\ge0$. Take its initial energy equal to
the energy of the restricted data. Inequality (2.7) there permits
nonnegative time tests supported up to zero and thus gives both forms.
Mass conservation cancels the affine term in the pressure potential;
the restricted initial energies increase to $E_{\delta,0}$.

Extend $(\rho_R,u_R)$ by zero for norms and for the continuity equation.
For each $T<\infty$, the energy inequality and Sobolev's inequality give,
uniformly in $R$,
\[
 \|m_R\|_{L^2(0,T;L^2)}
 \le \|\rho_R\|_{L^\infty(0,T;L^3)}\|u_R\|_{L^2(0,T;L^6)}
 \le C_{\delta,T}.
\]
The momentum equation is tested only on $B_R$. By
Lemma~\ref{lem:ball-pressure}, its pressure satisfies
\eqref{r4:fixed-pressure-moment} uniformly in $R$.

\emph{Passage to the whole space.}
For any compactly supported momentum test, its support lies in
$\R\times B_R$ for large $R$. Use the proof of
Theorem~\ref{cp:compactness} with fixed pressure $p_\delta$:
$m_R\in L^2L^2$ gives \eqref{cp:momentum_gain}, and
\eqref{r4:fixed-pressure-moment} gives local uniform integrability of
$p_\delta(\rho_R)$. Monotonicity of both pressure powers gives
\[
 \operatorname{Cov}(p_\delta(X),\vartheta_k(X))
 \ge\frac{\delta}{2}\,\mathbb E
         |\vartheta_k(X)-\mathbb E \vartheta_k(X)|^{\beta+1}.
\]
The renormalized continuity equation and the entropy argument yield
\[
 \begin{aligned}
 \rho_R&\to\rho &&\text{in }L^\beta_{\rm loc},\\
 p_\delta(\rho_R)&\to p_\delta(\rho)&&\text{in }L^1_{\rm loc},\\
 \rho_R|u_R|^2&\rightharpoonup\rho|u|^2&&\text{in }L^1_{\rm loc}.
 \end{aligned}
\]
The initial densities converge in $L^1\cap L^\beta$, and the initial
energies converge. Thus the trace and renormalization arguments in
Section~\ref{sec:compactness} apply with both powers of the density.

For the spatial tails, \eqref{r4:fixed-pressure-moment} provides
the pressure-moment bound required in
Proposition~\ref{r4:energy-tightness} with $\theta=1$, and
\[
 \int_I\int_{\{\rho_R>A\}}\rho_R|u_R|^2
 \le \int_I\|\rho_R\mathbf1_{\{\rho_R>A\}}\|_{3/2}\|u_R\|_6^2
 \le C_{\delta,I} A^{1-2\beta/3}.
\]
Propositions~\ref{r4:energy-tightness} and~\ref{r4:differential-energy}
give $E_{\delta,R}\to E_\delta$ in $L^1(I)$ and the differential
energy inequality; compare \cite{DellOroFeireisl2015}.
At time zero, use \eqref{cp:weighted_energy} and lower
semicontinuity. A diagonal choice of $R$ gives the global solution.
\end{proof}

\subsection{Approximation of the initial data}

Fix $\beta>\max\{4,\gamma\}$ and set
\begin{equation}\label{eq:approx-pressure}
 \delta_j=j^{-\beta-1},\qquad
 p_j(z)=az^\gamma+\delta_jz^\beta,\qquad
 \Theta_j(z)=\frac{a}{\gamma-1}z^\gamma+\frac{\delta_j}{\beta-1}z^\beta .
\end{equation}
With the usual zero value on vacuum, define
\begin{equation}\label{eq:prepared-data}
 \begin{gathered}
 w_0=\begin{cases}m_0/\sqrt{\rho_0},&\rho_0>0,\\0,&\rho_0=0,\end{cases}\\
 \rho_{0j}=\min\{\rho_0,j\}\,\mathbf1_{B_j},\qquad
 m_{0j}=\sqrt{\rho_{0j}}\,w_0 .
 \end{gathered}
\end{equation}
Finite initial energy gives $w_0\in L^2$. By dominated convergence
for the density and H\"older's inequality for the momentum,
\begin{equation}\label{eq:data-convergence}
 \begin{aligned}
 \rho_{0j}&\to\rho_0 &&\text{in }L^1\cap L^\gamma,\\
 m_{0j}&\to m_0 &&\text{in }L^1\cap L^{2\gamma/(\gamma+1)} .
 \end{aligned}
\end{equation}
Indeed, $\sqrt{\rho_{0j}}\to\sqrt{\rho_0}$ in $L^2\cap L^{2\gamma}$.
Moreover,
\begin{align}
 \int_{\R^3}\frac{|m_{0j}|^2}{\rho_{0j}}\dd x
 &=\int_{B_j\cap\{\rho_0>0\}}|w_0|^2\dd x
 \longrightarrow \int_{\R^3}\frac{|m_0|^2}{\rho_0}\dd x,
       \label{eq:prepared-kinetic}\\
 0\le\delta_j\int_{\R^3}\rho_{0j}^\beta\dd x
 &\le\delta_j j^{\beta-\gamma}\int_{\R^3}\rho_0^\gamma\dd x
 \longrightarrow0. \label{eq:prepared-artificial}
\end{align}
It follows that
\begin{equation}\label{eq:prepared-total}
 E_{0j}=\int_{\R^3}
 \left(\frac{|m_{0j}|^2}{2\rho_{0j}}+\Theta_j(\rho_{0j})\right)\dd x
 \longrightarrow E_0 .
\end{equation}
Only finiteness of $\|\rho_{0j}\|_\beta$ at each fixed $j$ is needed.
The data in \eqref{eq:prepared-data} require no further regularization.

\subsection{Uniform bounds on the whole space}

Apply Lemma~\ref{r4:fixed-pressure-energy} on $(0,j)$ with coefficient
$\delta_j$, data \eqref{eq:prepared-data}, and zero external force.
The resulting solutions satisfy the exact momentum equation and
renormalized continuity, with their initial terms. At fixed $j$,
finite mass and kinetic energy give $m_j\in L^\infty L^1$. A fixed universal radial mass
cutoff with $|\nabla\chi_R|\le C/R$ has error at most
$CR^{-1}\|m_j\|_{L^1((0,T)\times\R^3)}$. Thus
\begin{equation}\label{eq:approx-mass}
 \int_{\R^3}\rho_j(t,x)\,\dd x=\int_{\R^3}\rho_{0j}(x)\,\dd x .
\end{equation}
The usual pressure potential normalized at $z=1$ is
\begin{equation}\label{eq:affine-potential}
 z\int_1^z\frac{p_j(s)}{s^2}\,\mathrm ds
 =\Theta_j(z)-\left(\frac{a}{\gamma-1}+\frac{\delta_j}{\beta-1}\right)z .
\end{equation}
The linear term cancels by \eqref{eq:approx-mass}. For almost every $t$,
\begin{equation}\label{eq:approx-energy}
 \begin{aligned}
 &\int_{\R^3}\left(\frac12\rho_j|u_j|^2+\Theta_j(\rho_j)\right)(t)\,\dd x\\
 &\qquad+\int_0^t\int_{\R^3}
 \bigl(\mu|\nabla u_j|^2+(\mu+\lambda)|\diver u_j|^2\bigr)\,\dd x\,\dd s
 \le E_{0j}.
 \end{aligned}
\end{equation}
The solutions also satisfy the differential energy inequality. By
\eqref{eq:prepared-total}, for each fixed $T$ the following bounds
are independent of $j>T$:
\begin{equation}\label{eq:uniform-approx}
 \begin{gathered}
 \rho_j\in L^\infty_t(L^1_x\cap L^\gamma_x),\qquad
 u_j\in L^2_t\dot H^1_x,\\
 m_j,\ p_j(\rho_j),\ \rho_j|u_j|^2\in L^\infty_tL^1_x .
 \end{gathered}
\end{equation}
Indeed $p_j(z)\le\max\{\gamma-1,\beta-1\}\Theta_j(z)$.

\subsection{Renormalized continuity at the initial time}

For fixed $j$, the artificial pressure gives $\rho_j\in L^\infty L^\beta$
and hence $\rho_j\in L^\infty L^3$ by interpolation with $L^1$. Thus
\begin{equation}\label{eq:fixed-index-momentum}
 m_j\in L^2_tL^2_x,\qquad
 \rho_j u_j\otimes u_j\in L^1_tL^r_x,\qquad
 \frac1r=\frac1\beta+\frac13<1 .
\end{equation}
These bounds justify the test in Proposition~\ref{r4:global-pressure};
their constants need not be uniform in $j$.

At fixed $j$ extend $(\rho_j,u_j)$ to $t<0$ by $(\rho_{0j},0)$.
The weak initial trace gives the continuity equation across zero.
Since $\rho_j\in L^2_{\rm loc}$ and $\nabla u_j\in L^2_{\rm loc}$,
the spatial commutator gives the renormalized equation there. Its
initial term is $b(\rho_{0j})$.

\subsection{Limit as the artificial pressure vanishes}

\begin{proof}[Completion of the proof of Theorem~\ref{thm:nine-eighths}]
Fix $\gamma>9/8$ and $I'\Subset I\Subset(0,\infty)$.
The whole-space approximations satisfy the balance laws
\eqref{ao:balance} and the bounds \eqref{eq:uniform-approx}.
Consequently \eqref{ao:projection} holds uniformly in $j$.
Proposition~\ref{ao:gain} and Lemma~\ref{kt:tails} give
\[
 \begin{gathered}
 0<2\alpha<\min\left\{\tfrac16,\tfrac{8\gamma-9}{4\gamma+1}\right\},
 \qquad \kappa=\frac{(\gamma-1)\alpha}{1+2\alpha},\\
 \|I_{1-\alpha}m_j\|_{L^2(I'\times\R^3)}\le C_{I',I,\alpha},\qquad
 \int_{I'}\int_{\{\rho_j>A\}}\rho_j|u_j|^2
       \le C_{I',I,\alpha}A^{-\kappa}\qquad(A\ge2).
 \end{gathered}
\]
Since $(1+|\xi|^2)^{-1+\alpha}\le|\xi|^{-2+2\alpha}$,
\begin{equation}\label{eq:gain-for-compactness}
 \|m_j\|_{L^2(I';H^{-1+\alpha})}
 \le \|I_{1-\alpha}m_j\|_{L^2(I'\times\R^3)}\le C_{I',I,\alpha} .
\end{equation}
Multiplication by a fixed smooth compactly supported spatial cutoff gives
\eqref{cp:momentum_gain}. On $[0,L]$, the artificial pressure is
bounded by $\delta_j L^\beta\to0$. The other hypotheses of
Theorem~\ref{cp:compactness} are \eqref{eq:data-convergence},
\eqref{eq:prepared-total}, \eqref{eq:approx-energy} and
\eqref{eq:uniform-approx}. Hence, locally on interior cylinders,
\[
 p_j(\rho_j)\to a\rho^\gamma\quad\text{in }L^1,\qquad
 \rho_j u_j\otimes u_j\rightharpoonup\rho u\otimes u
       \quad\text{in }L^1 .
\]
The same theorem gives renormalized continuity with the prescribed
initial data, mass conservation, the traces
\eqref{eq:weak-traces} and \eqref{eq:strong-initial-density}, and
\eqref{cp:limit_energy}.

For positive times, \eqref{eq:gain-for-compactness} and
Proposition~\ref{r4:global-pressure} give
\eqref{r4:global-pressure-moment}. With the high-density bound
just obtained, Proposition~\ref{r4:energy-tightness} gives
\[
 E_j\to E\quad\text{in }L^1(I').
\]
Proposition~\ref{r4:differential-energy} now gives
\eqref{r4:energy-definition} and the initial energy trace.

A diagonal subsequence defines $(\rho,u)$ on $(0,\infty)\times\R^3$.
It satisfies \eqref{ao:balance} with $P=a\rho^\gamma$ and the mass,
energy and dissipation bounds \eqref{ao:natural}. Applying
Proposition~\ref{ao:gain} and Lemma~\ref{kt:tails} to this solution gives
\eqref{intro:signed-gain}--\eqref{intro:signed-tails}.
\end{proof}

\appendix

\section{Test functions for pressure estimates}\label{app:compact-tools}
\subsection{The effective viscous flux}\label{app:flux-test}

\begin{proof}[Proof of Lemma~\ref{cp:rough_test}]
Choose $\chi_0\in C_c^\infty(\R^3)$ equal to one near the spatial
supports of $f$ and $\zeta$. Set $v=\chi_0u$, so $fv=fu$ and
$v\in L^2H^1$. Fix once and for all a nonnegative standard spatial
mollifier $J_1\in C_c^\infty(B_1)$ with $\int J_1=1$, and write
$J_\varepsilon(x)=\varepsilon^{-3}J_1(x/\varepsilon)$ and
$f^\varepsilon=J_\varepsilon*f$. Constants $C_{J_1}$ below depend
only on this profile and the dimension. By
\eqref{cp:rough_cap_equation},
\[
 \partial_t \partial_i\Delta^{-1}f^\varepsilon
   =-R_{ik}J_\varepsilon*(fv_k)+\partial_i\Delta^{-1}(J_\varepsilon*g) .
\]
For fixed $\varepsilon$, this derivative belongs to
$L^2_tC^q_{\rm loc,x}$ for every integer $q\ge0$
(Lemma~\ref{cp:potential}). The spatial derivatives of the potential
are locally bounded. Time mollification on an interior interval
therefore gives admissible tests: the derivatives in time converge
in $L^2_tC^q_{\rm loc,x}$, the derivatives in space in
$L^s_tC^q_{\rm loc,x}$ for $s<\infty$. This suffices since
$m,\rho u\otimes u,P\in L^\infty_tL^1_x$ and $S\in L^2$.

Both $v_k R_{ik}f^\varepsilon$ and $R_{ik}(v_kf^\varepsilon)$ belong
to $L^2H^1$. By Lemma~\ref{cp:kinetic_duality},
\[
 \langle \zeta m_i,v_k R_{ik}f^\varepsilon\rangle
   =\int\zeta \rho u_i u_k R_{ik}f^\varepsilon .
\]
For this identity smooth $R_{ik}f^\varepsilon$ in time first; its
spatial derivatives are bounded, so dominated convergence applies
in $L^2H^1$ and against $\rho u\otimes u$. The time and convective terms now contain
\[
 \zeta m_i[v_k,R_{ik}]f^\varepsilon
 +\zeta m_i R_{ik}
       \bigl(v_kf^\varepsilon-J_\varepsilon*(v_kf)\bigr)
 +\zeta m\cdot \nabla\Delta^{-1}(J_\varepsilon*g).
\]
For $q_\varepsilon=(J_\varepsilon*f)\,v-J_\varepsilon*(vf)$,
the translation estimate for $H^1$ functions gives
\[
 \|q_\varepsilon\|_{L^2_{t,x}}
 \leq C_{J_1}\varepsilon\|f\|_\infty\|\nabla v\|_{L^2_{t,x}}.
\]
Differentiating the convolution gives
\[
 \partial_\ell q_\varepsilon(x)
 =(J_\varepsilon*f)(x)\,\partial_\ell v(x)
  +\int \partial_\ell J_\varepsilon(z)\cdot
       [v(x)-v(x-z)]f(x-z)\,\dd z .
\]
Thus
\[
 \|q_\varepsilon\|_{L^2H^1}\le C_{J_1}\|f\|_\infty\|v\|_{L^2H^1}.
\]
Interpolating with the $L^2$ error gives
\begin{equation}\label{cp:convolution_error}
 q_\varepsilon\longrightarrow0
 \quad\hbox{strongly in }L^2_tH^{1-\sigma}_x .
\end{equation}
Pair \eqref{cp:convolution_error} with $\zeta m$ after applying
$R_{ik}$. For $3/p<\sigma$ and $p<\infty$, one has
$f^\varepsilon(t)\to f(t)$ in $L^p$ a.e. in $t$, uniformly bounded
in that norm. Lemma~\ref{cp:commutator} and dominated convergence give
\[
 [v_k,R_{ik}]f^\varepsilon\longrightarrow[v_k,R_{ik}]f
 \quad\hbox{strongly in }L^2H^{1-3/p}.
\]
For the remaining terms, Lemma~\ref{cp:potential} and the potential
kernel ($p>3$) give
\[
 \begin{aligned}
 \nabla\Delta^{-1}(J_\varepsilon*g)&\to\nabla\Delta^{-1}g &&\text{in }L^2H^1,\\
 R_{ik}f^\varepsilon&\to R_{ik}f &&\text{in }L^2,\\
 \nabla\Delta^{-1}f^\varepsilon&\to\nabla\Delta^{-1}f
                   &&\text{in }L^2_tL^\infty_{\rm loc,x}.
 \end{aligned}
\]
Use bounded a.e. convergence against $P$ for $Pf^\varepsilon$.
The terms with derivatives of $\zeta$ converge by the
$L^\infty_tL^1_x$ bounds on $m,\rho u\otimes u,P$ and the $L^2$ bound on $S$.
\end{proof}

\subsection{A pressure estimate independent of the ball}
\label{app:bogovskii}
\begin{proof}[Proof of Lemma~\ref{lem:ball-pressure}]
Constants with a $\delta$ subscript have the dependencies specified in
Lemma~\ref{lem:ball-pressure}; the estimates below are uniform in $R$
and in the density truncation level $N$ unless stated otherwise.
On a fixed finite time interval, uniformly in $R\ge1$,
\[
 \begin{aligned}
 \rho_R&\in L^\infty(L^1\cap L^\beta),&
 u_R&\in L^2\dot H^1,\\
 m_R=\rho_Ru_R&\in L^\infty L^{6/5}\cap L^2L^2,&
 \rho_R u_R\otimes u_R&\in L^1L^{3/2}.
 \end{aligned}
\]
These bounds follow from energy, interpolation and
$\|u_R\|_6\le C\|\nabla u_R\|_2$. The momentum equation is used only on $B_R$.

Let $\mathcal B_R$ be the Bogovski\u{\i} operator on $B_R$, and
$\langle f\rangle_R$ the average of $f$ on the ball. The field
$\mathcal B_R(f-\langle f\rangle_R)$ has zero trace and divergence
$f-\langle f\rangle_R$. We use
\[
 \|\nabla\mathcal B_R(f-\langle f\rangle_R)\|_q\le C_q\|f\|_q
 \quad(1<q<\infty),\qquad
 \|\mathcal B_R\diver v\|_2\le C\|v\|_2 .
\]
For the second estimate interpret $\diver v$ on $H^1(B_R)/\mathbb R$ by
$\langle\diver v,\phi\rangle=-\int v\cdot\nabla\phi$. Approximate
$v$ in $L^2$ by compactly supported smooth fields. The estimate extends
$\mathcal B_R\diver$ by continuity for this zero-flux interpretation. If the zero extension of $v$ is
divergence free, shrinking it inside the ball and then mollifying
approximates it by divergence-free test fields in $L^2$. Thus the
extension depends only on this divergence functional. Both constants are independent of $R$ by scaling;
see \cite[Lemma~4]{DellOroFeireisl2015}.

Use the truncations $\vartheta_N$ of \eqref{cp:density_caps}. In the following
calculation, the only source in the renormalized equation is
\[
 g_N=\bigl(\vartheta_N(\rho_R)-\rho_R \vartheta_N'(\rho_R)\bigr)\diver u_R,\qquad
 0\le \vartheta_N(\rho_R)-\rho_R \vartheta_N'(\rho_R)\le\rho_R .
\]
At fixed $N$, $\vartheta_N$ is bounded and its flux has zero normal trace.
Indeed, set $d_R(x)=\operatorname{dist}(x,\partial B_R)$.
For a boundary cutoff of thickness $\varepsilon$, the flux error is
at most
\[
 C_N\int_0^T\int_{\{d_R<\varepsilon\}}
              \frac{|u_R|}{d_R}\,\dd x\,\dd t\longrightarrow0
\]
by Hardy's inequality and Cauchy--Schwarz. Thus we may test
renormalized continuity up to the boundary.

For $\zeta\in C_c^\infty(0,T)$ with $0\le\zeta\le1$, test momentum with
\[
 \zeta\,\mathcal B_R
       \bigl[\vartheta_N(\rho_R)-\langle \vartheta_N(\rho_R)\rangle_R\bigr].
\]
Renormalized continuity gives the derivative of the spatial factor:
\[
 -\mathcal B_R\diver(\vartheta_N(\rho_R)u_R)
       +\mathcal B_R\bigl[g_N-\langle g_N\rangle_R\bigr].
\]
The following bounds are uniform in $N,R$:
\[
 \begin{aligned}
 \|\mathcal B_R\diver(\vartheta_N(\rho_R)u_R)\|_{L^2L^2}
 &\le C\|m_R\|_{L^2L^2},\\
 \|\mathcal B_R(g_N-\langle g_N\rangle_R)\|_{L^2L^2}
 &\le C\|g_N\|_{L^2L^{6/5}}
 \le C\|\rho_R\|_{L^\infty L^3}\|\diver u_R\|_{L^2L^2}.
 \end{aligned}
\]
The second line uses $W^{1,6/5}_0(B_R)\hookrightarrow L^2(B_R)$.
For $\Psi=\mathcal B_R[\vartheta_N(\rho_R)-\langle \vartheta_N(\rho_R)\rangle_R]$,
the remaining pairings are bounded by
\[
 \begin{aligned}
 \|m_R\|_{L^2L^2}\|\partial_t\Psi\|_{L^2L^2}
 &+\|\rho_Ru_R\otimes u_R\|_{L^1L^{3/2}}\|\nabla\Psi\|_{L^\infty L^3}\\
 &+\|S(u_R)\|_{L^2L^2}\|\nabla\Psi\|_{L^2L^2}
 +C_\zeta\|m_R\|_{L^\infty L^{6/5}}\|\Psi\|_{L^\infty L^6}.
 \end{aligned}
\]
Here $\nabla\Psi\in L^\infty(L^2\cap L^3)$ and
$\Psi\in L^\infty L^6$, uniformly by the Bogovski\u{\i} and Sobolev estimates.

For admissibility at fixed $R,N$, extend $\zeta\Psi$ by zero in space
and time. Its spatial trace vanishes, so its distributional divergence
is the zero extension of $\zeta\cdot(\vartheta_N(\rho_R)-\langle \vartheta_N(\rho_R)\rangle_R)$.
For $0<\eta<1$ put
\[
 \Phi_\eta(t,x)=\eta\cdot(\zeta\Psi)(t,x/\eta),
 \qquad \supp_x\Phi_\eta\subset\overline B_{\eta R},\qquad
 \|\diver\Phi_\eta\|_\infty\le 2\|\zeta\|_\infty\|\vartheta_N\|_\infty .
\]
Mollify in space and time with radius less than $(1-\eta)R/2$ and
less than the distance to the time endpoints. First remove the mollifier,
then let $\eta\uparrow1$. These smooth, compactly supported tests
converge to $\zeta\Psi$ with time derivatives in $L^2L^2$ and gradients
in $L^2L^2$. Dilation and convolution also give strong convergence of the
gradients in $L^3_{t,x}$, with a uniform $L^\infty_tL^3_x$ bound.
Along a subsequence they converge in $L^3_x$ at a.e. time, which suffices
against the convective flux in $L^1_tL^{3/2}_x$. The divergences converge a.e. along a
subsequence with the displayed uniform bound, so the pressure term
converges against $p_\delta(\rho_R)\in L^1$.

The term containing the spatial mean of $\vartheta_N(\rho_R)$ is bounded by
\[
 \frac{\|\rho_R\|_{L^\infty L^1}}{|B_R|}
       \int_0^T\int_{B_R} \zeta\,p_\delta(\rho_R)\le C_{\delta,T}
       \qquad(R\ge1).
\]
Take $\zeta=1$ on $I\Subset(0,T)$. Since $\vartheta_N(z)\uparrow z$, the
claimed estimate, with constant $C_{\delta,I,T}$, follows by monotone
convergence.
\end{proof}

\section{The projected determinant estimate}
\label{app:middle-tools}
We prove the following quantitative form of
Proposition~\ref{rt:optimized-CI}. Write
\(\log_+s=\max\{0,\log s\}\). We fix universal cutoff profiles for
the localizations below. Dependence on the mollifier profile $J_1$ is
displayed explicitly. Under the proposition's hypotheses, for \(d,h>0\),
\begin{equation}\label{rt:free-D}
 \begin{aligned}
 \mathcal G_h(\mathcal A)\le C_{J_1}\bigg[&
 \|\mathcal F\|_1^{3/2}
 +d^{1/2}\left(\int_{\R^4}\tr\mathcal A\right)^{3/2}\\
 &+d^{-1/2}\|\mathcal E\|_2^{3/2}
       \left(1+\log_+\frac{d^{-1/3}}h\right)^{3/4}\bigg].
 \end{aligned}
\end{equation}
The constant is independent of $d,h$ and of the supports.

\subsection{The periodic Monge--Amp\`ere equation}
\label{app:projection-details}

On \(\mathbb T_L^3=\R^3/(L\mathbb Z^3)\) we use Lebesgue measure.
For a mean-zero function,
\[
 \|f\|_{\dot H^s(\mathbb T_L^3)}^2
 =L^3\sum_{n\ne0}\left|\frac{2\pi n}{L}\right|^{2s}|\widehat f(n)|^2.
\]
For a matrix, sum this expression over the components.

Let \(A_h\) be smooth, symmetric and nonnegative on \(\mathbb T_L^3\).
For \(\varepsilon>0\), let
\[
 f_\varepsilon=\det(A_h+\varepsilon\mathrm{Id}_3)^{1/2}.
\]
Take \(L^{-3}\le d/2\). We solve
\begin{equation}\label{rt:MA-D}
 \begin{gathered}
 \det(\mathrm{Id}_3+\nabla^2\psi)
 =1-\frac{L^{-3}}d+\frac{f_\varepsilon}{d\int_{\mathbb T_L^3}f_\varepsilon},\\
 \mathrm{Id}_3+\nabla^2\psi>0,\qquad
 \int_{\mathbb T_L^3}\psi=0 .
 \end{gathered}
\end{equation}
The right-hand side is smooth, is at least \(1/2\), and has mean one.
Hence a smooth periodic solution exists by
\cite[Theorem~2.2]{Li1990}, after rescaling.

\emph{Bounds for the gradient.}
We claim that
\begin{equation}\label{rt:MA-D-basic}
 \Delta\psi\ge-Cd^{-1}L^{-3},\qquad
 \|\Delta\psi\|_1\le Cd^{-1},\qquad
 \|\nabla\psi\|_\infty\le Cd^{-1/3}.
\end{equation}
Indeed, the arithmetic--geometric mean inequality gives
\[
 3+\Delta\psi\ge 3(1-d^{-1}L^{-3})^{1/3}.
\]
Since \(\int\Delta\psi=0\), the \(L^1\) bound follows.
The gradient of the periodic Green function is bounded by
\(C(\operatorname{dist}(y,0)^{-2}+L^{-2})\). Thus weak Young's
inequality gives
\begin{equation}\label{app:green-weak}
 \|\nabla\psi\|_{L^{3/2,\infty}(\mathbb T_L^3)}
 \le C\|\Delta\psi\|_1 .
\end{equation}

Write \(\ell_\psi=\|\nabla\psi\|_\infty\) in this paragraph, and suppose \(\ell_\psi>0\).
On each coordinate circle, \(\partial_i\psi\) has mean zero and
\(\partial_{ii}\psi\ge-1\). Integrating along the circle gives
\(\|\partial_i\psi\|_\infty\le L/2\), so \(\ell_\psi\le\sqrt3 L/2\).
Choose \(y_0\) with \(|\nabla\psi(y_0)|=\ell_\psi\). Convexity of
\(|y|^2/2+\psi\) on \(\R^3\) gives
\[
 \nabla\psi(y_0+z)\cdot z\ge\nabla\psi(y_0)\cdot z-|z|^2.
\]
It follows that \(|\nabla\psi(y_0+z)|\ge \ell_\psi/2\) on the cone
\[
 0<|z|<\ell_\psi/4,\qquad \nabla\psi(y_0)\cdot z>3\ell_\psi|z|/4 .
\]
The cone fits in one coordinate chart and has volume \(c\ell_\psi^3\).
Its contribution to the weak \(L^{3/2}\) norm is at least
\(c \ell_\psi\cdot(\ell_\psi^3)^{2/3}=c \ell_\psi^3\). Now \eqref{app:green-weak} gives
the last inequality in \eqref{rt:MA-D-basic}.

\emph{Bounds for the mollified test.}
For every ball of radius \(s<L/8\), the divergence theorem and
\eqref{rt:MA-D-basic} give
\begin{equation}\label{app:trace-local}
 \int_{B_s}|\Delta\psi|
 =\int_{\partial B_s}\nabla\psi\cdot n+2\int_{B_s}(\Delta\psi)_-
 \le C\bigl(d^{-1/3}s^2+d^{-1}L^{-3}s^3\bigr).
\end{equation}
Set \(r_0=d^{-1/3}\) and take \(L>16(r_0+h)\). If \(h<r_0\), then
\eqref{app:trace-local} and the global \(L^1\) bound give
\[
 \int_{B_s}|J_h*\Delta\psi|\le C_{J_1}
 \begin{cases}
 d^{-1/3}h^{-1}s^3+d^{-1}L^{-3}s^3,&s<h,\\
 d^{-1/3}s^2+d^{-1}L^{-3}s^3,&h\le s<r_0,\\
 d^{-1},&s\ge r_0 .
 \end{cases}
\]
The kernel of \((-\Delta)^{-1/2}\) on mean-zero periodic functions is
bounded in absolute value by \(C(\operatorname{dist}(y,0)^{-2}+L^{-2})\).
Both periodic kernel bounds follow by integrating the heat kernel
after subtracting its constant mode. Integrating the displayed ball
bound against \(2s^{-3}\dd s\), first up to \(h\), then to \(r_0\), and
then beyond \(r_0\), bounds the kernel potential of
\(|J_h*\Delta\psi|\) by
\[
 C_{J_1}d^{-1/3}\bigl(1+\log(r_0/h)\bigr)+Cd^{-1}L^{-2}.
\]
Multiplication by \(\|J_h*\Delta\psi\|_1\le Cd^{-1}\) now bounds the
squared \(\dot H^{-1/2}\) norm. If \(h\ge r_0\), use instead
\(\|J_h*\Delta\psi\|_\infty\le C_{J_1}d^{-1}h^{-3}\) for radii below \(h\)
and the \(L^1\) bound above \(h\). The squared norm is then at most
\[
 C_{J_1}d^{-2}(h^{-2}+L^{-2})\le C_{J_1}d^{-4/3}.
\]
Finally, \(\sum_{i,j}|\xi_i\xi_j|^2=|\xi|^4\). Thus
\begin{equation}\label{rt:trace-D}
 \|J_h*\nabla^2\psi\|_{\dot H^{-1/2}(\mathbb T_L^3)}
 =\|J_h*\Delta\psi\|_{\dot H^{-1/2}(\mathbb T_L^3)}
 \le C_{J_1}d^{-2/3}\sqrt{1+\log_+(d^{-1/3}/h)}.
\end{equation}

For the terms involving a spatial cutoff, we also need
\begin{equation}\label{rt:W-D}
 \|J_h*\nabla\psi\|_{\dot H^{-1/2}(\mathbb T_L^3)}
 \le C_{J_1}d^{-1}\sqrt{1+\log(L/h)}.
\end{equation}
Indeed \(|\widehat{\Delta\psi}(n)|\le Cd^{-1}L^{-3}\).
At \(n\ne0\), the multiplier taking \(\Delta\psi\) to
\(\nabla\psi\) has modulus \(L/(2\pi|n|)\). Hence the squared norm
in \eqref{rt:W-D} is bounded by
\[
 Cd^{-2}\sum_{n\ne0}|n|^{-3}
       |\widehat{J_1}(2\pi hn/L)|^2
 \le C_{J_1}d^{-2}\bigl(1+\log(L/h)\bigr).
\]
These constants do not depend on \(\varepsilon\) or \(A_h\).

\subsection{Proof of the line integral estimate}

\begin{lemma}[Localization on a periodic cell]\label{app:periodic-extension}
Let \(f\) be smooth, \(L\)-periodic and of mean zero. If
\(\chi\in C_c^\infty((-1/2,1/2)^3)\), then
\begin{equation}\label{app:extension-bound}
 \|\chi(\cdot/L)f\|_{\dot H^{-1/2}(\R^3)}
 \le C_\chi\|f\|_{\dot H^{-1/2}(\mathbb T_L^3)} .
\end{equation}
The constant is independent of \(L\).
\end{lemma}
\begin{proof}
It suffices to take \(L=1\): both norms scale by \(L^2\).
Multiplication by \(\chi\) followed by periodization is bounded
from \(H^s(\R^3)\) to \(H^s(\mathbb T^3)\) for \(s=0,1\), hence for
\(s=1/2\) by interpolation. Duality bounds the \(H^{-1/2}(\R^3)\) norm
of \(\chi f\). This controls the homogeneous norm for \(|\xi|\ge1\).
For the remaining frequencies,
\[
 |\widehat{\chi f}(\xi)|
 \le \|f\|_{H^{-1/2}(\mathbb T^3)}
      \|\chi(y)e^{-\mathrm i y\cdot\xi}\|_{H^{1/2}(\mathbb T^3)}
 \le C_\chi\|f\|_{\dot H^{-1/2}(\mathbb T^3)},\qquad |\xi|<1 .
\]
Here \(f\) has mean zero, and \(|\xi|^{-1}\) is integrable on \(B_1\).
\end{proof}

\begin{proof}[Proof of Proposition~\ref{rt:optimized-CI}]
\emph{Localization.}
Choose smooth cutoffs \(0\le\chi_R\uparrow1\) on \(\R^4\), with
\(\chi_R=1\) on \(B_R\), support in \(B_{2R}\), and
\(|\nabla\chi_R|\le C/R\). The product rule gives
\begin{equation}\label{app:cutoff-balance}
 \operatorname{Div}\bigl(\chi_R\cdot(\mathcal A-\mathcal E)\bigr)
 =\chi_R\mathcal F+\mathcal A\nabla\chi_R-\mathcal E\nabla\chi_R .
\end{equation}
For any compactly supported matrix field $V$, distributional
integration by parts gives
\[
 \diver_{\omega^\perp}\bigl(\Pi_\omega X_\omega V\Pi_\omega^{\top}\bigr)
 =\Pi_\omega X_\omega(\Div V).
\]
Indeed, lift a transverse test field to one constant along the lines;
compact support of $V$ permits a cutoff equal to one near its support.
Use this identity with $V=\chi_R\cdot(\mathcal A-\mathcal E)$.
Convolve with $J_h$ and periodize on a cell of side $L=16R$, with
$R$ sufficiently large for the preceding estimates. In \eqref{rt:MA-D} use
\[
 A_h=J_h*\bigl[\Pi_\omega X_\omega(\chi_R\mathcal A)
                         \Pi_\omega^{\top}\bigr].
\]

\emph{Testing the equation on \(\omega^\perp\).}
The matrix arithmetic--geometric mean inequality and
\eqref{rt:MA-D} give
\[
 (A_h+\varepsilon\mathrm{Id}_3):(\mathrm{Id}_3+\nabla^2\psi)
 \ge 3d^{-1/3}\left(\int_{\mathbb T_L^3}f_\varepsilon\right)^{-1/3} f_\varepsilon.
\]
Consequently
\begin{equation}\label{rt:D-left}
 3d^{-1/3}\left(\int_{\mathbb T_L^3}f_\varepsilon\right)^{2/3}
 \le \int_{\mathbb T_L^3}
     (A_h+\varepsilon\mathrm{Id}_3):(\mathrm{Id}_3+\nabla^2\psi).
\end{equation}
The trace term is at most \(\int_{\R^4}\tr\mathcal A+3\varepsilon L^3\), and
\(\int\varepsilon\Delta\psi=0\). For the remaining term,
test the projected equation with \(\nabla\psi\). Since \(J_h\) is even,
we may put the convolution on the test. The source terms
\(\chi_R\mathcal F\) and \(\mathcal A\nabla\chi_R\) are bounded in
\(L^1\) by \(\|\mathcal F\|_1\) and
\(CR^{-1}\int_{\R^4}\tr\mathcal A\). For the two terms containing
\(\mathcal E\), pair in \(\dot H^{1/2}\) and \(\dot H^{-1/2}\).
The projected fields are supported in $B_{2R}\subset\omega^\perp$.
Thus one may multiply each periodic test by $\chi(\cdot/L)$ with
$\chi=1$ on $B_{1/8}$ before taking Euclidean Sobolev duality.
Lemma~\ref{app:periodic-extension} then gives the periodic norms
in \eqref{rt:trace-D} and \eqref{rt:W-D}.
Together with \eqref{rt:MA-D-basic}, this gives
\begin{align}
 \left(\int_{\mathbb T_L^3}f_\varepsilon\right)^{2/3}
 \le C_{J_1}\bigg[&
 d^{1/3}\left(\int_{\R^4}\tr\mathcal A+3\varepsilon L^3\right)
       +\|\mathcal F\|_1+R^{-1}\int_{\R^4}\tr\mathcal A \notag\\
 &+d^{-1/3}\sqrt{1+\log_+(d^{-1/3}/h)}
       \,\|X_\omega(\chi_R\mathcal E)\|_{\dot H^{1/2}(\omega^\perp)}
       \notag\\
 &+d^{-2/3}\sqrt{1+\log(L/h)}
       \,\|X_\omega(\mathcal E\nabla\chi_R)\|_{\dot H^{1/2}(\omega^\perp)}
 \bigg].
 \label{rt:directional-D}
\end{align}

\emph{Integration over directions.}
For \(f\in C_c^\infty(\R^4)\), the Fourier transform of \(X_\omega f\)
on \(\omega^\perp\) is the restriction of \(\widehat f\) to this plane.
Thus
\begin{align}
 \int_{\mathbb S^3}\|X_\omega f\|_{\dot H^{1/2}(\omega^\perp)}^2\dd\omega
 &=C\int_{\R^4}|\xi|\,|\widehat f(\xi)|^2
       \left(\int_{\mathbb S^3}\delta_0(\omega\cdot\xi)\dd\omega\right)\dd\xi
       \notag\\
 &=C\|f\|_2^2.
 \label{app:slice}
\end{align}
Here \(\delta_0\) is the one-dimensional Dirac mass, and the integral
in parentheses is \(C/|\xi|\). By density the identity holds for
compactly supported \(L^2\) fields, component by component. In
particular the two Sobolev norms in \eqref{rt:directional-D} need
only be finite for almost every $\omega$; no bound uniform in
$\omega$ is used. Orthogonal projection does not increase the component norm.

We first keep $R,d,h$ fixed and remove the positive regularization;
only afterwards do we let $R\to\infty$, and then choose $d$ as a
function of $h$. In particular, at fixed \(R,d,h\), let \(\varepsilon\downarrow0\) in
\eqref{rt:directional-D}. Raise both sides to the power \(3/2\) and
integrate in \(\omega\). H\"older's inequality on \(\mathbb S^3\) and
\eqref{app:slice} bound the two stress norms in angular \(L^2\) by
\(C\|\mathcal E\|_2\) and \(C\|\mathcal E\|_2/R\), respectively. In particular, the last term in
\eqref{rt:directional-D} has angular \(L^2\) norm at most
\begin{equation}\label{app:cutoff-stress-error}
 C_{J_1}d^{-2/3}\frac{\|\mathcal E\|_2}{R}\sqrt{1+\log(L/h)}
 \longrightarrow0\qquad(R\to\infty;\ d,h\text{ fixed}).
\end{equation}
The term \(R^{-1}\int_{\R^4}\tr\mathcal A\) also tends to zero. On the left, the projected
mollifications of \(\chi_R\mathcal A\) increase in matrix order:
\(\mathcal A\ge0\) and \(\chi_R\uparrow1\). Their determinants increase
as well. Monotone convergence gives \eqref{rt:free-D}.

For $0<h\le1$, take $d=(1+|\log h|)^{3/4}$. Since $d^{-1/3}\le1$,
\eqref{rt:free-D} gives
\[
 \begin{aligned}
 \mathcal G_h(\mathcal A)\le C_{J_1}\bigg[&\|\mathcal F\|_1^{3/2}\\
 &+\left(\left(\int_{\R^4}\tr\mathcal A\right)^{3/2}
          +\|\mathcal E\|_2^{3/2}\right)
          (1+|\log h|)^{3/8}\bigg],
 \end{aligned}
\]
as asserted in Proposition~\ref{rt:optimized-CI}.
\end{proof}

\section{Counting the time intervals in the momentum cover}
\label{app:momentum-cover}

\begin{proof}[Proof of Lemma~\ref{ao:cover-lemma}]
We use the constant convention of Section~\ref{ao:section}.
The spatial convolution of \eqref{ao:balance} and the kernel bound
\[
 |\nabla\eta_s|\le Cs^{-1}\eta_s
\]
give
\[
 |\partial_tm_s|
 \le \frac Cr\,\eta_s*(|\rho u\otimes u+P\mathrm{Id}_3|+|S|),
 \qquad r/2<s<2r.
\]
For fixed \(r>0\) the right side is bounded in \(L^1(I;L^\infty)\)
uniformly in \(s\in[r/2,2r]\). Hence \(f_j\in W^{1,1}(I)\) and
\begin{equation}\label{ao:variation}
 |f_j'(t)|
 \le \frac Cr\int w_j(x)\cdot
                \bigl(|\rho u\otimes u+P\mathrm{Id}_3|+|S|\bigr)(t,x)\,\mathrm dx .
\end{equation}
For every finite subcollection of intervals, bounded overlap and
\eqref{ao:envelope-mass} yield
\[
 \sum_j w_j\le C,\qquad \int w_j\le Cr^3,\qquad
 \Lambda r^3\ell\le C\sum_j\int_I f_j\le C_*.
\]
This proves \eqref{ao:weight-bounds}. On each interior interval
the variation of \(f_j\) is at least \(c_0\Lambda r^3/2\). An interval meeting
an endpoint contributes to the corresponding trace; by Fatou's lemma,
\[
 \sum_j\bigl(f_j(t_-)+f_j(t_+)\bigr)\le C_* .
\]
Summing \eqref{ao:variation} and the endpoint contributions gives
\[
 rN\le \frac {C_*}{\Lambda r^3}
 \left(\iint w\cdot\bigl(|\rho u\otimes u+P\mathrm{Id}_3|+|S|\bigr)+r\right).
\]
Adding $\ell$ gives \eqref{ao:raw-count}. By \eqref{ao:natural} and Cauchy--Schwarz,
\[
 \iint w|\rho u\otimes u+P\mathrm{Id}_3|\le C_*,\qquad
 \iint w|S|\le C\|S\|_2 (r^3\ell)^{1/2}\le C_* \Lambda^{-1/2}.
\]
Together these inequalities give \eqref{ao:coarse-complexity} for \(\Lambda\ge1\).
For each fixed \(\Lambda>0\) they bound the number of intervals independently
of the finite subcollection. The full collection is therefore finite.
\end{proof}

\end{document}